\documentclass[11pt]{amsart}
\usepackage{amssymb,amscd}
\numberwithin{equation}{section}

\newtheorem{theorem}{Theorem}[section]

\newtheorem{corollary}[theorem]{Corollary}
\newtheorem{proposition}[theorem]{Proposition}
\newtheorem{lemma}[theorem]{Lemma}
\theoremstyle{definition}
\newtheorem{definition}[theorem]{Definition}
\theoremstyle{remark}
\newtheorem{remark}[theorem]{Remark}

\newcommand\cS{\mathcal{S}}
\newcommand\cC{\mathcal{C}}

\newcommand{\W}{\mathcal{W}}

\newcommand{\R}{\mathbb{R}}
\newcommand{\C}{\mathbb{C}}
\newcommand{\Z}{\mathbb{Z}}

\newcommand{\bbT}{\mathbb{T}}
\newcommand\lie[1]{\mathfrak{#1}}

\newcommand{\fh}{\lie{h}}
\newcommand{\fg}{\lie{g}}

\newcommand{\fk}{\lie{k}}
\newcommand{\ft}{\lie{t}}

\def    \inv    {^{-1}}

\newcommand{\SP}        {\operatorname{Sp}}
\newcommand{\Iso}	{\operatorname{Iso}}

\newcommand{\cL}{\underline{\Z_G}}

\newcommand{\labell}[1]	{\label{#1}}

\begin{document}

\title{Contact 	toric manifolds }
\author{Eugene Lerman}
\address{Department of
Mathematics, University of Illinois, Urbana, IL 61801}
\email{lerman@math.uiuc.edu}
\date{\today}
\thanks{Partially supported by NSF grant DMS - 980305,
the American Institute of Mathematics and R. Kantorovitz.}

\begin{abstract} 
We complete the classification of compact connected contact toric
manifolds initiated by Banyaga and Molino and by Galicki and Boyer. As an
application we prove the conjectures of Toth and Zelditch on toric
integrable systems on the $n$-torus and the 2-sphere.
\end{abstract}

\maketitle
\tableofcontents

\section{Introduction}

The main goal of this paper is to complete the classification of
compact co-oriented contact manifolds with completely integrable
contact torus actions, that is, of compact contact toric manifolds.
See Theorem~\ref{main-theorem} below for a precise description of our
classification.  The study of completely integrable systems on contact
manifolds was initiated by Banyaga and Molino \cite{BM1, BM2, B} who,
in particular, gave a partial classification of contact toric
manifolds.  Boyer and Galicki gave a description of contact toric
manifolds satisfying the assumption that the image of the moment map
lands in an open half-space \cite{BG}.

Our motivation for studying contact toric manifolds comes from the
recent work of Toth and Zelditch \cite{TZ} who introduced the notion
of a toric integrable geodesic flow.  For every toric integrable
geodesic flow on a manifold $M$ the co-sphere bundle $S(T^*M)$ has the
structure of a contact toric manifold.  Thus a classification of
contact toric manifolds  is useful for understanding of toric
integrable geodesic flows.  For example in
\cite{LS} it was shown that if a geodesic flow on the standard torus 
$\bbT^n = \R^n /\Z^n$ is toric integrable then the corresponding metric
must be flat thereby proving a conjecture in \cite{TZ}.  The argument
in \cite{LS} uses a topological classification of contact toric
manifolds with a non-free torus action.

Alternatively, since there is a one-to-one correspondence between
symplectic cones and contact manifolds, our classification is a
classification of symplectic toric manifolds that are symplectic cones
over a compact base.  Compact symplectic toric manifolds are well
understood thanks to the work of Delzant \cite{D}.  All of such
manifolds turn out to be smooth projective toric varieties.  In
general compact toric manifolds have been studied extensively from
many different points of view.  It seems to us that the class of
non-compact symplectic toric manifolds we are describing should have
interesting applications beyond integrable geodesic flows.

A classification of contact toric manifolds up to equivariant
contactomorphisms allows us in this paper to positively answer two
 questions posed in \cite{TZ}:

\begin{theorem}\labell{thm1.1}
Up to  appropriate isomorphisms the only toric integrable action of
the $n$-torus on the co-sphere bundle of the $n$-torus is the standard
action, that is,  the lift of multiplication on the left to the co-sphere
bundle.
\end{theorem}

\begin{theorem}\labell{thm1.2}
Up to  appropriate isomorphisms the only toric integrable action of
the $2$-torus on the co-sphere bundle of the 2-sphere is the standard
action, that is the action generated by the geodesic flow of the round
metric and by the rotations about a fixed axis.
\end{theorem}

We strongly suspect that there are preciously few manifolds that admit
toric integrable geodesic flows.  In fact the only examples we know
are tori and the low-dimensional spheres, and it would be surprising
if there are others.

Note that it is an entirely different problem to classify all toric
integrable {\bf metrics} on a given manifold. As was mentioned earlier
for any toric integrable metric on a torus has to be flat and
conversely any flat metric is toric integrable.  But even for the 2-sphere
the question is open.


\subsection*{A note on notation}  

If $U$ is a subspace of a vector space $V$ we denote its annihilator
in the dual vector space $V^*$ by $U^\circ$.  Thus $U^\circ = \{ \ell
\in V^* \mid \ell |_U = 0\}$.

Throughout the paper the Lie algebra of a Lie group denoted by a
capital Roman letter is denoted by the same small letter in the
fraktur font: thus $\fg$ denotes the Lie algebra of a Lie group $G$
etc.  The vector space dual to $\fg$ is denoted by $\fg^*$. The
identity element of a Lie group is denoted by 1.  The natural pairing
between $\fg$ and $\fg^*$ is denoted by $\langle \cdot, \cdot
\rangle$.

When a Lie group $G$ acts on a manifold $M$ we denote the action by an
element $g\in G$ on a point $x\in M$ by $g\cdot x$; $G\cdot x$ denotes
the $G$-orbit of $x$ and so on.  The vector field induced on $M$ by an
element $X$ of the Lie algebra $\fg$ of $G$ is denoted by $X_M$ and
the diffeomorphism induced by $g\in G$ on $M$ by $g_M$.  Thus in this
notation $g\cdot x = g_M (x)$.  The isotropy group of a point $x\in M$
is denoted by $G_x$; the Lie algebra of $G_x$ is denoted by $\fg_x$
and is referred to as the isotropy Lie algebra of $x$.  Recall that
$\fg_x = \{ X \in \fg\mid X_M (x) = 0\}$.

If a Lie group $G$ is a torus we denote the integral lattice of $G$ by
$\Z_G$ and the dual weight lattice by $\Z_G^*$.  Recall that $\Z_G =
\ker (\exp :\fg \to G)$, $\Z_G^* = \text{Hom}_\Z (\Z_G, \Z)$.

If $X$ is a vector field and $\tau$ is a tensor, then $L_X \tau$
denotes the Lie derivative of $\tau$ with respect to $X$.

If $P$ is a principal $G$-bundle then $[p, m]$ denotes the point in the
associated bundle $P\times _G M = (P\times M)/G$ which is the orbit of
$(p,m) \in P\times M$.

The symbol $\alpha$ always denotes a contact form and $\xi$ always
denotes a co-oriented contact structure.

\subsection*{Acknowledgments} This work was written during a visit of  the
 American Institute of Mathematics and of the University of California
 at Berkeley.  The author is grateful to the University of California
 at Berkeley for providing access to its computing facilities.

\section{Group actions on contact manifolds}
In this section we recall the notion of a co-oriented contact structure
and of a contact form, and discuss the definition of a moment map for
group actions on contact manifolds.

\begin{definition}
Recall that a 1-form $\alpha$ on a manifold $M$ is {\bf contact} if
$\alpha _x \not = 0$ for any $x\in M$, so $\xi = \ker \alpha$ is a
codimension-1 distribution, and if additionally $d\alpha |_\xi$ is
non-degenerate. Thus the vector bundle $\xi \to M$ necessarily has
even-dimensional fibers, and the manifold $M$ is necessarily
odd-dimensional.
\end{definition}

\begin{definition}
A codimension-1 distribution $\zeta$ on a manifold $M$ is {\bf
co-orientable} if its annihilator $\zeta^\circ \subset T^*M$ is an
oriented line bundle, i.e., has a nowhere vanishing global section. It
is {\bf co-oriented} if one component $\zeta^\circ _+$ of $\zeta^\circ
\smallsetminus 0$ ($\zeta^\circ$ minus the zero section) is chosen.
\end{definition}

\begin{definition}
A {\bf co-oriented contact structure} $\xi$ on a manifold $M$ is a
co-oriented codimension-1 distribution such that $\xi^\circ
\smallsetminus 0$  is a symplectic submanifold of the cotangent bundle 
$T^*M$ (the cotangent bundle is given the canonical symplectic form).
We denote the chosen component of $\xi^\circ \smallsetminus 0$ by
$\xi^\circ _+$ and refer to it as the {\bf symplectization} of $(M,
\xi)$.
\end{definition}

\begin{remark}
It is a standard fact that $\xi \subset TM$ is a co-oriented contact
structure if and only if there is a contact form $\alpha$ with $\ker
\alpha = \xi$ (given $\xi$ choose $\alpha$ to be a section of 
$\xi ^\circ \smallsetminus 0 \to M$).  If $f$ is any
function on $M$ then $e^f \alpha$ defines the same contact structure
$\xi$ and conversely, if $\alpha$ and $\alpha'$ are two contact forms
defining the same contact sturcture with the same co-orientation, then
$\alpha' = e^f \alpha$ for some function $f$.  That is, a co-oriented
contact structure is the same thing as a conformal class of contact
forms.
\end{remark}

\begin{remark}
In this paper $\alpha$ always denotes a contact form and $\xi$ always
denotes a co-oriented contact structure (with a co-orientation
understood).  We will refer to a pair $(M, \alpha)$ or to a pair $(M,
\xi)$ as a (co-oriented) {\bf contact manifold}.
\end{remark}

\begin{lemma}\labell{lemma2.6}
Suppose a Lie group $G$ acts properly on a manifold $M$ preserving a
co-oriented codimension-1 distribution $\zeta$ and its co-orientation.
That is, suppose the lifted action of $G$ on $T^*M$ preserves a
component $\zeta^\circ_+$ of $\zeta^\circ \smallsetminus 0$.  Then
there is a $G$-invariant 1-form $\beta$ on $M$ such that $\zeta = \ker
\beta$ and $\beta (M) \subset\zeta^\circ_+$. 
\end{lemma}
\begin{proof}  
To choose $\beta$ first choose any 1-form $\tilde{\beta}$ on $M$ with
$\ker \tilde{\beta} = \zeta$ and $\tilde{\beta} (M)\subset \zeta^\circ
_+$.  If $G$ is compact, we then let $\beta$ to be the average $\int
_G (g_M)^* \tilde{\beta} \, dg$ of $\tilde{\beta}$ over the group $G$.
If $G$ is not compact, then due to the existence of slices we may
assume that $M = G\times _K V$ for some representation $K \to GL(V)$
of a compact Lie group $K$.  We  average the restriction
$\tilde{\beta}|_V : V \to T^* (G\times _K V)$ over $K$ and then extend
it to all of $M$ by $G$-invariance.
\end{proof}

Therefore given a proper action of a Lie group $G$ on a contact
manifold $(M, \xi = \ker \alpha)$ which preserves the co-orientation,
we may (and do) assume that the contact form $\alpha$ is $G$-invariant.

\begin{definition}
If a Lie group $G$ acts on a manifold $M$ preserving a 1-form
$\beta$, the corresponding $\beta$-{\bf moment map } $\Psi_\beta :M
\to \fg^*$ determined by $\beta$ is defined by
\begin{equation}\labell{stupid-eq}
\langle \Psi _\beta (x) , X \rangle = \beta _x (X_M (x))
\end{equation}
for all $x\in M$ and all vectors $X$ in the Lie algebra $\fg$ of $G$,
where $X_M$ denotes the vector field induced by $X$: $X_M (x) =
\frac{d}{dt} |_{t=0} (\exp tX)\cdot x$.
\end{definition}
 If $d\beta$ is a symplectic form then, up to a sign
convention, $\Psi_\beta$ is a symplectic moment map.  If $\alpha$ a
contact form then $\Psi_\alpha$ is a candidate for a contact moment
map.  Note however that if $f$ is a $G$-invariant function, then $e^f
\alpha$ is also a contact form defining the same contact distribution,
while clearly $\Psi_{e^f \alpha} = e^f \Psi _\alpha$. That is, this
definition of the moment map depends on a particular choice of a
contact form and not just on the contact structure.

Fortunately there is also a notion of a contact moment map that
doesn't have this problem.  Namely, suppose again that a Lie group $G$
acts on a manifold $M$ preserving a co-oriented contact structure
$\xi$.  The lift of the action of $G$ to the cotangent bundle then
preserves a component $\xi^\circ_+$ of $\xi^\circ \smallsetminus 0$.
The restriction $\Psi =\Phi|_{\xi^\circ_+}$ of the moment map $\Phi$
for the action of $G$ on $T^*M$ to $\xi^\circ _+$ depends only on the
action of the group and on the contact structure.  Moreover, since
$\Phi: T^*M \to \fg^*$ is given by the formula 
$$
\langle \Phi (q, p), X\rangle = \langle p, X_M (q)\rangle 
$$ 
for all $q\in M$, $p\in T^*_q M$ and  $X \in \fg$, we see
that if $\alpha$ is any invariant contact form with $\ker \alpha =
\xi$ and $\alpha (M) \subset \xi^\circ_+$ 
then $\langle \alpha^* \Psi (q, p), X\rangle = \langle \alpha ^*
\Phi (q, p), X\rangle = \langle \alpha _q, X_M (q) \rangle = \langle
\Psi _\alpha (q), X\rangle$.    Thus $\Psi \circ \alpha = \Psi _\alpha$, 
that is, $\Psi =  \Phi|_{\xi^\circ_+}$ is a ``universal'' moment map.

There is another reason why the universal moment map $\Psi :\xi^\circ_+
\to \fg^*$ is a more natural notion of the moment map than the one
given by (\ref{stupid-eq}).  The vector fields induced by the action
of $G$ preserving the contact distribution $\xi$ are contact.  The
space of contact vector fields is isomorphic to the space of sections
of the bundle $TM/\xi \to M$ (a choice of a contact form identifies
$TM/\xi$ with $M\times \R$ and contact vector fields with functions).
Thus a contact group action gives rise to a linear map
\begin{equation}\labell{eq**}
\fg \to \Gamma (TM/\xi), \quad X\mapsto X_M \mod \xi .
\end{equation}
The moment map should be the transpose of the map (\ref{eq**}).  The
total space of the bundle $(TM/\xi)^*$ naturally maps into the space
dual to the space of sections $\Gamma (TM/\xi)$: 
$$ 
(TM/\xi)^* \ni
\eta \mapsto \left(s \mapsto \langle \eta, s (\pi (\eta))
\rangle\right), 
$$ 
where $\pi : (TM/\xi)^* \to M$ is the projection and $\langle \cdot,
\cdot \rangle$ is the paring between the corresponding fibers of
$(TM/\xi)^*$ and $TM/\xi$.  In other words, the transpose $\Psi:
(TM/\xi)^* \to \fg^*$ of (\ref{eq**}) should be given by
\begin{equation}\labell{eq***}
\langle \Psi (\eta), X\rangle 
	= \langle \eta, X_M (\pi (\eta))\!  \mod \xi \rangle
\end{equation}
Under the identification $\xi^\circ \simeq (TM/\xi)^*$, the equation
above becomes 
$$
\langle \Psi (q, p), X\rangle = \langle p, X_M (q)\rangle
$$ 
for all $q\in M$, $p\in \xi_q^\circ$ and $X\in \fg$, which exactly the
definition of $\Psi $ given earlier as the restriction of the moment
map for the lifted action of $G$ on the cotangent bundle $T^*M$.
Thus part of the above discussion can be summarized as 

\begin{proposition}
Let $(M, \xi)$ be a co-oriented contact manifold with an action of a
Lie group $G$ preserving the contact distribution and its
co-orientation.  Suppose there exists an invariant 1-form $\alpha$
with $\ker \alpha = \xi$ and $\alpha (M) \subset \xi^\circ_+$ (c.f.\
Lemma~\ref{lemma2.6}).  Then the $\alpha$-moment map $\Psi_\alpha $ for the
action of $G$ on $(M, \alpha)$ and the moment map $\Psi$ for the
action of $G$ on the symplectization $\xi^\circ _+$ are related by
$$
\Psi \circ \alpha = \Psi _\alpha.
$$ 
Here $\xi^\circ _+$ is the component of $\xi \smallsetminus 0$ containing 
the image of $\alpha :M \to \xi ^\circ$.
\end{proposition}

\begin{remark}
We will refer to $\Psi: \xi^\circ_+ \to \fg^*$ as the moment map for
the action of a Lie group $G$ on a co-oriented contact manifold $(M,
\xi = \ker \alpha)$, that is, as the {\bf contact moment map}.
It is easy to show that $\Psi$ is $G$-equivariant with respect to the
given action of $G$ on $M$ and the coadjoint action of $G$ on $\fg^*$.
Hence for any invariant contact form $\alpha$ the corresponding
$\alpha$-moment map $\Psi_\alpha :M \to \fg^*$ is also
$G$-equivariant.
\end{remark}

Later in the paper we will need a version of contact reduction due to
Albert \cite{Albert} and, independently, to Geiges \cite{Geiges}:
\begin{lemma}[Contact quotients]\labell{lemma-contact-quotient}
  Suppose a Lie group $G$ acts on a manifold $M$ preserving a 1-form
  $\beta$.  Let $\Psi _\beta : M \to \fg^* $ denote the corresponding
  moment map.  Suppose $\Psi_\beta \inv (0) $ is a manifold and
  suppose that $G$ acts freely and properly on $\Psi_\beta \inv (0)$.
  Then $\beta $ descends to a 1-form $\beta _0$ on $M_0 := \Psi_\beta
  \inv (0)/G$.

If $\beta$ is contact than $\beta_0$ is contact as well.  Moreover the
manifold $M_0$ and the contact structure on $M_0$ defined by $\beta_0$
depends only on the contact structure defined by $\beta$ and not on the form
$\beta$ itself.
\end{lemma}

\begin{proof}
It's easy to see that $\beta |_{\Psi_\beta\inv (0)} $ is basic and
hence descends.  For a proof that $\beta_0$ is contact if $\beta $ is
contact see \cite{Albert} or \cite{Geiges}.  Note that the zero level
set $\Psi_\beta\inv (0)$ is the set of points of $M$ where $G$-orbits
are tangent to the contact distribution $\ker \beta$.  Thus $M_0$
depends only on $\ker \beta$.  It is not hard to see that $\ker \beta
_0$ depends only on $\ker \beta$ as well.
\end{proof}

We now begin the study of contact toric manifolds.
\begin{definition}
An action of a torus $G$ on a contact manifold $(M, \xi)$ is {\bf
completely integrable} if it is effective, preserves the contact
structure $\xi$ and if $2\dim G = \dim M +1$.

A {\bf contact toric $G$-manifold} is a co-oriented contact manifold
$(M, \xi)$ with a completely integrable action of a torus $G$ (but see
Remark~\ref{remark-below} below).
\end{definition}

Note that if an action of a torus $G$ on $(M, \xi)$ is completely
integrable, then the action of $G$ on a component $\xi^\circ_+$ of
$\xi ^\circ \smallsetminus 0$ is a completely integrable Hamiltonian
action and thus $\xi_+^\circ$ is a symplectic toric manifold (for more
information on symplectic toric manifolds and orbifolds see \cite{D}
and \cite{LT}).

\begin{lemma}\labell{lemma-nonzero}
Suppose an action of a torus $G$ on a contact manifold $(M, \xi)$ is
completely integrable.  Then zero is not in the image of the contact
moment map $\Psi : \xi^\circ _+ \to \fg^*$.
\end{lemma}
\begin{proof}
Suppose not. Then for some point $x\in M$ the orbit $G\cdot x$ is
tangent to the contact distribution $\xi$.  Therefore the tangent space
$\zeta _x := T_x (G\cdot x)$ is isotropic in the symplectic vector
space $(\xi _x , \omega_x)$ where $\omega_x = d\alpha _x|_\xi$ and
$\alpha$ is a $G$-invariant contact form with $\ker \alpha = \xi$ and
$\alpha (M) \subset \xi^\circ_+$.

We now argue that this forces the action of $G$ not to be effective.
More precisely we argue that the slice representation of the connected
component of identity $H$ of the isotropy group $G_x$ of the point $x$ is
not effective. The group $H$ acts on $\xi_x$ preserving the symplectic
form $\omega_x$ and preserving $\zeta_x$.  Since
$\zeta_x$ is isotropic, $\xi_x = (\zeta_x^\omega/ \zeta_x) \oplus
(\zeta_x \times \zeta_x ^*)$ as a symplectic representation of $H$.
Here $\zeta_x ^\omega$ denotes the symplectic perpendicular to
$\zeta_x$ in $(\xi_x, \omega_x)$.  Note that since $G$ is a torus, the
action of $H$ on $ \zeta_x$ (and hence on $\zeta_x^*$) is trivial.

Observe next that the dimension of the symplectic vector space $V :=
\zeta_x^\omega/ \zeta_x$ is $\dim \xi_x - 2 \dim \zeta_x = 
\dim M -1 - 2 (\dim G - \dim H) = (\dim M - 1) - (\dim M + 1) + 2 \dim H
 = 2 \dim H - 2$.  
On the other hand, since $H$ is a compact connected Abelian group acting
symplecticly on $V$, its image in the group of symplectic linear
transformations $\SP (V)$ lies in a maximal torus $T$ of a maximal
compact subgroup of $\SP (V)$.  Since the maximal compact subgroup of 
$\SP (V)$ is the unitary group $U(n)$, $n= \dim V/2$, the dimension of 
$T$ is $n =
\dim H -1$.  Therefore the representation of $H$ on $V$ (and hence of $G_x$)
 is not faithful.  Since the fiber at $x$ of the normal bundle of
 $G\cdot x$ in $M$ is $T_xM/\xi _x \oplus \xi_x/\zeta_x \simeq \R
 \oplus V \oplus
\zeta_x^*$, the slice representation of $G_x$ is not faithful.
Consequently the action of $G$ in not faithful in a neighborhood of an
orbit $G\cdot x$.  Contradiction.
\end{proof}

Suppose $(M, \xi)$ is a contact toric $G$-manifold.  Fix an inner
product on the Lie algebra $\fg$ of $G$ and thereby on the dual space
$\fg^*$.  There exists then a unique $G$-invariant contact form
$\alpha$ defining $\xi$ and its co-orientation and furthermore
satisfying the normalizing condition that $||\Psi_\alpha (x) || = 1$
at all $x\in M$ (if $\alpha'$ is any contact form defining $\xi$ let
$\alpha_x = \frac{1}{||\Psi_{\alpha'} (x)||} \alpha' _x$ for $x\in M$;
by Lemma~\ref{lemma-nonzero} this makes sense).

\begin{remark}\labell{remark-below}
From now on we fix an inner product on a torus $G$.  This allows us to
normalize contact forms on contact toric $G$-manifolds as above.
Also, since the moment map $\Psi_\alpha: M \to \fg^*$ contains all the
information about the action of $G$ on $(M, \xi = \ker \alpha)$ we can
and will think of a {\bf contact toric $G$-manifold} as a triple $(M,
\alpha, \Psi_\alpha)$ where $\alpha $ is normalized so that 
$|| \Psi_\alpha (x)|| = 1 $ for all $x\in M$.
\end{remark}

\begin{definition}
Let $(M, \xi)$ be a co-oriented contact manifold with an action of a
Lie group $G$ preserving the contact structure $\xi$ and its
co-orientation.  Let $\Psi: \xi^\circ_+ \to \fg^*$ denote the
corresponding moment map.  We define the {\bf moment cone} $C(\Psi)$
to be the set $$ C(\Psi) : = \Psi (\xi^\circ_+)\cup \{0\}.  $$ Note
that if $\alpha $ is a $G$-invariant contact form with $\xi = \ker
\alpha$ and $\alpha (M) \subset \xi^\circ_+$,  then 
$$
C(\Psi ) = \{t f \mid f \in \Psi _\alpha (M), \, t \in [0, \infty)\},
$$
where  $\Psi_\alpha :M\to \fg^*$ denote the $\alpha$-moment map.
\end{definition}

\begin{definition} \label{def2.15}
Two contact toric $G$-manifolds 
$(M, \alpha, \Psi_\alpha)$ and $(M, \alpha', \Psi_{\alpha'})$ are {\bf
isomorphic } if there exists a $G$-equivariant co-orientation
preserving contactomorphism $\varphi: M \to M'$.  We will refer to
such a map $\varphi$ as an {\bf isomorphism} between $(M, \alpha,
\Psi_\alpha)$ and $(M, \alpha', \Psi_{\alpha'})$.

We denote the group of isomorphisms of $(M, \alpha, \Psi_\alpha)$ by 
$\Iso (M, \alpha, \Psi_\alpha) = \Iso (M)$
\end{definition}

\begin{remark}\label{rmrk2.15}
Note that if $\varphi: M \to M'$ is an isomorphism of two contact
toric $G$-manifolds $(M, \alpha, \Psi_\alpha)$ and $(M, \alpha',
\Psi_{\alpha'})$, then $\varphi^* \alpha' = e^f \alpha$ for some
$G$-invariant function $f\in C^\infty (M)$.  Consequently $\varphi^*
\Psi_{\alpha'} = e^f \Psi _\alpha$.  But for any $x\in M$, $1 = ||
\Psi_{\alpha'} (\varphi (x))|| = ||e^{f(x)} \Psi _\alpha (x)|| =
e^{f(x)}1 = e^{f(x)}$.  Thus $\varphi^*\alpha ' = \alpha$ and
$\varphi^* \Psi_{\alpha'} = \Psi_\alpha$.
\end{remark}

\begin{definition}[Good cones] \label{def-good-cones}
 Let $\fg^*$ be the dual of the Lie algebra 
of a torus $G$.  Recall that a subset $C\subset \fg^* $ is a {\bf
rational polyhedral cone} if there exists a finite set of vectors
$\{v_i\}$ in the integral lattice $\Z_G$ of $G$ such that 
$$ 
C =
\bigcap \{ \eta \in \fg^* \mid \langle \eta, v_i\rangle \geq 0\}. 
$$
Of course it is no loss of generality to assume that the set $\{v_i\}$
is {\bf minimal}, i.e., that for any index $j$ 
$$ 
C \not = \bigcap
_{i\not =j} \{ \eta \in \fg^* \mid \langle \eta, v_i\rangle \geq 0\},
$$ 
and that each vector $v_i$ is {\bf primitive}, i.e., $sv_i \not \in
\Z_G$ for $s \in (0, 1)$.  Therefore we make these two assumptions.

A rational polyhedral cone $C = \bigcap_{i=1}^N \{ \eta \in \fg^* \mid
\langle \eta, v_i\rangle \geq 0\}$, $\{v_i\} \subset \Z_G$ with
non-empty interior is {\bf good}\footnote{for the purposes of this
paper} if the annihilator of a linear span of a codimension $k$ space,
$0<k < \dim G$ is the Lie algebra of a subtorus $H$ of $G$ and the
normals to the face form a basis of the integral lattice $\Z_H$ of
$H$.  That is, if 
$$
\{0\} \not = C \cap 
\bigcap _{j= 1}^k \{ \eta \in \fg^* \mid \langle \eta, v_{i_j}\rangle \geq 0\}
$$ 
is a face of $C$ for some $\{i_1, \ldots, i_k\} \subset \{1, \ldots, N\}$ then
\begin{equation}
\{ \sum _{j=1}^k a_j v_{i_j} \mid a_j \in \R\} \cap \Z_G = \{ \sum _{j=1} ^k m_j v_{i_j} \mid m_j \in \Z\}
\end{equation}
and $\{v_{i_j}\}$ is independent over $\Z$.
\end{definition}

We can now state the main result of the paper, a classification theorem.

\begin{theorem}\labell{main-theorem}
Compact connected contact toric (c.c.c.t.) $G$-manifolds $(M, \alpha,
\Psi_\alpha : M\to \fg^*)$ are classified as follows.\\

1. \ Suppose $\dim M = 3$ and the action of $G= \bbT^2$ is free. Then
$M$ is a principal $G$-bundle over $S^1$, hence is diffeomorphic to
$\bbT^3 = S^1\times \bbT^2$.  Moreover the contact form $\alpha$ is
$\cos nt \, d\theta_1 + \sin nt \, d \theta_2$ ($(t, \theta_1,
\theta_2) \in S^1\times \bbT^2$) for some positive integer $n$.\\

2. \ Suppose $\dim M = 3$ and the action of $G = \bbT^2$ is not free.
Then $M$ is diffeomorphic to a lens space (this includes $S^1 \times
S^2$) and, as a c.c.c.t. $G$-manifold, $(M, \alpha, \Psi_\alpha)$ is
classified by two rational numbers $r, q$ with $0\leq r < 1$, $r< q$.\\

3. \ Suppose $\dim M > 3$ and the action of $G$ is free.  Then $M$ is
a principal $G$-bundle over a sphere $S^d$, $d= \dim G-1$.  Moreover
each principal $G$-bundle over $S^d$ has a unique $G$-invariant
contact structure making it a c.c.c.t. $G$-manifold.\\

4. \ Suppose $\dim M > 3$ and the action of $G$ is not free.  Then the
moment cone of $(M, \alpha, \Psi_\alpha)$ is a good cone (cf.\
Definition~\ref{def-good-cones} above). Conversely, given a good cone
$C\subset \fg^*$ there is a unique c.c.c.t. $G$-manifold $(M, \alpha,
\Psi_\alpha)$ with moment cone $C$.
\end{theorem}

\begin{remark}
Since principal $n$-torus bundles over a manifold are in one-to-one
correspondence with second cohomology classes of the manifold with
coefficients in $\Z^n$ and since $H^2 (S^d, \Z^n) = 0$ for $d\not =
2$, the only interesting case of the part 3 of the theorem occurs when
$\dim G = 3$.  In this case the theorem asserts that each principal
$\bbT^3$ bundle over the 2-sphere (there are $\Z^3$ of them
altogether) carries a unique $\bbT^3$-invariant contact structure.
These contact manifolds were first constructed by Lutz \cite{Lutz}.  Their
symplectizations were also explicitly constructed by Bates \cite{Bates}.
\end{remark}

We will prove Theorem~\ref{main-theorem} over the course of the next
four sections.

\section{Local structure of contact toric manifolds}

The motivation for the following definition comes from
Lemma~\ref{lem.pre-iso.orb} below.  Note that given an embedding
$\iota :N \to M$ we do not distinguish between the vector bundles over
$N$ and over $\iota (N)$.

\begin{definition}
Let $(M, \xi = \ker \alpha)$ be a co-oriented contact manifold.  An
embedded submanifold $N\hookrightarrow M$ is {\bf pre-isotropic} if 
\begin{enumerate}
\item $N$ is transverse to the contact distribution $\xi$ and 
\item the distribution $\zeta  = TN \cap \xi$ is isotropic in the 
conformal symplectic vector bundle $(\xi, [\omega])$ where $[\omega]$
is the conformal class of $ d\alpha |_{\xi}$.
\end{enumerate}
\end{definition}

\begin{remark} Note that $\zeta = \ker (\iota ^* \alpha)$.  Note also that if
$\alpha' = e^f \alpha$, $f\in C^\infty (M)$, is another contact form
defining the contact structure $\xi$ then $\iota ^* \alpha' =
e^{\iota^* f}\iota^* \alpha$ and $d\alpha' |_\xi = e^f (d\alpha
|_\xi)$. 
\end{remark}

\begin{definition} 
Let $\iota: N\hookrightarrow (M, \xi = \ker \alpha)$ be a pre-isotropic
embedding.  We define the {\bf characteristic distribution} of the
embedding $\iota$ to be the co-oriented distribution $\zeta = TN \cap
\xi$. Equivalently we can think of $\zeta$ as the conformal class
$[\iota^*\alpha]$ of 1-forms.  We define the {\bf conformal symplectic
normal bundle } $(E, [\omega_E])$ of the embedding by $E = \zeta
^\omega /\zeta$ where $ \zeta ^\omega$ is the symplectic perpendicular
to $\zeta$ in the conformal symplectic vector bundle $(\xi, [\omega] =
[d\alpha |_\xi])$ and $[\omega_E]$ is the conformal class of
symplectic structures induced on $E$ by $[\omega]$.
\end{definition}

\begin{remark}
Suppose $\iota: N\hookrightarrow (M, \xi = \ker \alpha)$ is a
pre-isotropic embedding.  Suppose further that a Lie group $G$ acts
on $N$ and $M$ preserving the contact form $\alpha$ and making the
embedding $\iota$ equivariant.  Then $G$ preserves the
characteristic distribution $\zeta$ and acts on the conformal
symplectic normal bundle $E$ preserving the symplectic structure
$\omega_E$ and its conformal class.
\end{remark}
\begin{theorem}[Uniqueness of pre-isotropic embeddings]\labell{thm-unique-pre}
A pre-isotropic embedding is uniquely determined by its characteristic
distribution and its conformal symplectic normal bundle.

More specifically suppose $(M_j,\xi_j = \ker \alpha_j)$, $j=1,2$ are
two contact manifolds and $\iota_j :N \to (M_j,\xi_j)$, $j=1,2$ are two
pre-isotropic embeddings such that 
$$
\iota_1 ^*\alpha _1 = e^f \iota_2^* \alpha _2
$$
and
$$
(E_1, \omega_1) \simeq (E_2, e^h \omega_2) \quad 
\text{as symplectic vector bundles},
$$ where $f, h \in C^\infty (N)$ are two functions and $(E_1,
[\omega_1])$ and $(E_2, [\omega_2])$ are the conformal symplectic
normal bundles of the embeddings.

Then there exist neighborhoods $U_j$ of $\iota_j (N)$ in $M_j$
($j=1,2$) and a diffeomorphism $\varphi : U_1 \to U_2$ such that
$\iota_2 = \varphi \circ \iota_1$ and $\varphi^*\alpha_2 = e^g \alpha
_1$ for some $g\in C^\infty (U_1)$.

Moreover if a Lie group $G$ acts properly on $N$, $M_1$, $M_2$ making
the embeddings $\iota_j$ $G$-equivariant and if the action preserves
the contact structures, then we may choose the neighborhoods $U_1$,
$U_2$ to be $G$-invariant and the map $\varphi$ to be $G$-equivariant.
\end{theorem}

The proof of Theorem~\ref{thm-unique-pre} relies on the
following observation.

\begin{theorem}[Equivariant relative Darboux theorem] \labell{lemma2}
Let $N\hookrightarrow M$ be an embedded closed submanifold. Suppose
there exist on $M$ two contact structures $\xi^0 =\ker \alpha^0$ and
$\xi^1 = \ker \alpha ^1$ ($\alpha^0, \alpha^1$ are 1-forms) and a
function $f\in C^\infty (M)$ such that
$$
\alpha^0_x = e^{f(x)} \alpha^1_x \quad \text{for all } x\in N,
$$
and suppose that the 2-form
$$
\left( (1-t) d\alpha ^0_x + t d \alpha ^1 _x\right) |_{\xi^0_x}
$$
is nondegenerate for all $x\in N$ and all $t\in [0, 1]$.

Then there exist neighborhoods $U_0, U_1$ of $N$ in $M$ and a
diffeomorphism $\varphi : U_0 \to U_1$ such that $\varphi |_N = id_N$
and $\varphi ^*\alpha ^1 = e^h \alpha ^1$ for some $h \in C^\infty (U_0)$.

Moreover, if a Lie group $G$ acts properly on $M$ preserving $N$ and
the two contact forms $\alpha^0$, $\alpha^1$, then we can choose the
neighborhoods $U_0, U_1$ to be $G$-invariant and arrange for the map
$\varphi$ above to be $G$-equivariant.
\end{theorem}   

\begin{proof}
Consider the family of $G$-invariant 1-forms $\alpha ^t =t\alpha ^1 +
(1-t) \alpha ^0$, $t\in [0, 1]$.  For all $x\in N$ and all $t\in [0,
1]$ we have $\ker \alpha^t_x = \xi^0_x = \xi^1 _x$ and $d\alpha^t_x |_{\ker
\alpha^t_x}$ is nondegenerate.  Therefore the forms $\alpha^t$ are contact
in a neighborhood of $N$ for all $t$.  It is no loss of generality to
assume that this neighborhood is all of $M$.

Denote the Reeb vector field of $\alpha^t$ by $Y_t$.  Since the Reeb
vector field is uniquely defined by $\alpha ^t (Y_t) = 1$, $\iota
(Y_t) d\alpha ^t = 0$ and since $\alpha ^t$ is $G$-invariant, $Y_t$ is
$G$-invariant as well.

Define a time dependent vector field $X_t$ tangent to the contact
distribution $\xi^t = \ker \alpha^t$ by
$$
X_t = \left(d\alpha^t|_{\xi^t }\right) \inv (-\dot {\alpha}^t|_{\xi^t})
$$
where $\dot{\alpha} ^t = \frac{d}{dt} \alpha^t$ is the derivative with
respect to $t$.  Clearly $X_t$ is $G$-invariant.  Note that $X_t (x) =
0$ for all $x\in N$.  This is because $-\dot {\alpha}^t_x |_{\xi^t_x}
= (\alpha^0_x - \alpha ^1_x) |_{\xi^t_x} = (e^{f(x)} - 1 )\alpha^1_x
|_{\xi^1_x} =0$ for $x\in N$.
We claim that the Lie derivative of $\alpha^t$ with respect to $X_t$ 
satisfies 
\begin{equation}\labell{eq-Lie}
L_{X_t} \alpha^t = \dot{\alpha}^t (Y_t) \alpha ^t - \dot{\alpha}^t .
\end{equation}
Indeed, since $\alpha ^t(X_t) = 0$, $L_{X_t} \alpha^t = \iota (X_t)
d\alpha^t$.  By definition of $X_t$, $(\iota (X_t) d\alpha^t)|_{\xi^t}
= -\dot{\alpha }^t|_{\xi^t} = \left( \dot{\alpha}^t (Y_t) \alpha ^t -
\dot{\alpha}^t\right)|_{\xi^t}$.  On the other hand, 
$(\iota (X_t) d\alpha^t) (Y_t) = 0  = \dot{\alpha}^t (Y_t)\,
 1 -  \dot{\alpha}^t (Y_t) =  \dot{\alpha}^t (Y_t) \alpha ^t (Y_t) -
\dot{\alpha}^t (Y_t)$.  This proves (\ref{eq-Lie}).  Hence
\begin{equation} \labell{next-eq}
L_{X_t} \alpha^t + \dot{\alpha}^t= \dot{\alpha}^t(Y_t) \alpha ^t .
\end{equation}
Denote the isotopy generated by $X_t$ by $\varphi_t$. Since $X_t$
vanishes on $N$, $\varphi_t $ is defined for all $t\in [0, 1]$ on a
neighborhood of $N$. Since $X_t$ is $G$-invariant, the isotopy is
$G$-equivariant.  Let $g_t = \varphi_t^* (\dot{\alpha}^t(Y_t))$.  Then
$\frac{d}{dt}(\varphi_t^* \alpha^t) = \varphi_t^* (L_{X_t}\alpha^t +
\dot{ \alpha}^t) = g_t (\varphi^*_t
\alpha^t)$ by equation (\ref{next-eq}).  Therefore $\varphi^*_t \alpha
^t = \left(e^{\int_0^t g_s \, ds}\right) \varphi ^*_0 \alpha^0$.
In particular,
$$
\varphi^*_1 \alpha^1 = e^h \alpha^0,
$$
where $h = \int_0^1 g_s \, ds \in C^\infty (M)$.  Note also that since $X_t
$ is zero at the points of $N$, $\varphi_t$ fixes $N$ pointwise.
\end{proof}

\begin{proof}[Proof of Theorem~\ref{thm-unique-pre}]

Suppose $\iota :N \to (M, \xi = \ker \alpha)$ is a $G$-equivariant
pre-isotropic embedding. Choose a $G$-invariant almost complex
structure $J$ on $\xi$ compatible with $\omega = d\alpha |_{\xi}$.
This gives us a $G$-invariant inner product on the vector bundle $\xi
\to M$.  Extend it to a $G$-invariant Riemannian metric $g$ on $M$ by, say, 
declaring the Reeb vector field $Y$ of the invariant contact form
$\alpha$ to be of unit length and orthogonal to $\xi$ (as remarked
previously, since $\alpha$ is $G$-invariant and the Reeb vector field
$Y$ is uniquely determined by $\alpha$, $Y$ is $G$-invariant).

By construction of $J$, $\zeta := \ker (\alpha |_N)$ and $J\zeta$ are
$g$-perpendicular and $\zeta \oplus J\zeta$ is a symplectic subbundle
of $(\xi, \omega)$.  Let $E$ be the $\omega$-perpendicular to
$\zeta\oplus J\zeta$ in $\xi$.  The bundle $E$ is also $g$-perpendicular
to $\zeta\oplus J\zeta$; $E$ is isomorphic to the conformal symplectic
normal bundle of the embedding $\iota$.  We therefore have a
$G$-equivariant direct sum decomposition 
$$
\xi|_N =  \zeta \oplus J\zeta \oplus E =
 (TN \cap \xi) \oplus J\zeta \oplus E.
$$
Since $N$ is transverse to $\xi$  it follows that 
$$
TM|_N = TN \oplus J\zeta \oplus E.
$$

Note that the Reeb vector field $Y$ need not be tangent to $N$, and so
$J\zeta \oplus E$ need not be $g$-perpendicular to $TN$.  Never the
less, $J\zeta \oplus E$ is a topological normal bundle for the
embedding $\iota :N \to M$.  Therefore the restriction of the
$g$-exponential map $\exp : TM \to M$  to $J\zeta \oplus
E$ gives an open $G$-equivariant embedding of a neighborhood of the
zero section $N \hookrightarrow J\zeta \oplus E$ into $M$; the
embedding is identity on $N$.

Since $\zeta$ is isotropic in $(\xi, \omega)$, the map $J\zeta \to
\zeta^*$ defined by $v \mapsto \omega (v, \cdot)|_\zeta$ is an isomorphism. 
By composing the inverse of this map with $\exp$ we get a
$G$-equivariant map $\psi : \zeta^* \oplus E \to M$, which has the
following properties.  The map $\psi$ is the identity on $N$.  It is
an open embedding on a sufficiently small neighborhood of $N$ in
$\zeta^* \oplus E$.  For any point $(x, 0, 0)$ in the zero section of
$\zeta^* \oplus E$ we have
$$
\ker (\psi^* \alpha)_{(x,0,0)} = \zeta_x\oplus\zeta_x^* \oplus E_x
$$
and
$$
d (\psi^* \alpha)_{(x,0,0)}|_{\ker (\psi^* \alpha)_{(x,0,0)}} =
\omega_{\zeta_x\oplus\zeta_x^*} \oplus (\omega_E)_x, 
$$ 
where $\omega_{\zeta_x\oplus\zeta_x^*}$ denotes the canonical
symplectic form on $\zeta_x\oplus\zeta_x^*$ and $(\omega_E)_x =
d\alpha _x |_{E_x}$.

Now consider two pre-isotropic embeddings $\iota_j : N \to (M_j,
\zeta_j)$, $j=1,2$ satisfying the hypotheses of the theorem.  
Since $\iota_1^*\alpha_1 = e^f \iota_2^* \alpha _2$, $\ker \iota_1
^*\alpha _1 = \ker \iota_2^* \alpha _2$.  Denote this distribution by
$\zeta$.  Let $\sigma : E_1 \to E_2$ denote the vector bundle
isomorphism with $\sigma ^*\omega_2 = e^h \omega_1$.  The map $\tau =
id \oplus \sigma: \zeta^* \oplus E_1 \to \zeta^* \oplus E_2$ has the
property that at the points $(x,0,0)$ of the zero section $N$
\begin{equation}
d\tau_{(x, 0, 0)} =
 id_{T_x N}\oplus id_{\zeta_x^*} \oplus \sigma |_{(E_1)_x}.
\end{equation}

 Consider the two  $G$-equivariant maps $\psi_j: \zeta_j^* \oplus E_j
\to M$, $j=1,2$ given by the construction at the beginning of the proof.
We have, for any point $x\in N$,
$$
\aligned
\ker ((\psi_2\circ \tau)^*\alpha_2)_{(x,0,0)} & 
= (d\tau _{(x, 0, 0)})\inv (\zeta_x \oplus \zeta^*_x \oplus (E_2)_x)  \cr
&= \zeta_x \oplus \zeta^*_x \oplus (E_1)_x\cr
&= \ker (\psi_1^* \alpha_1)_{(x,0,0)}.
\endaligned
$$
Also
$$
d (\psi_1^* \alpha_1)_{(x,0,0)}|_{\ker (\psi_1^* \alpha_1)_{(x,0,0)}} =
\omega_{\zeta_x \oplus \zeta_x^*} \oplus (\omega_1)_x
$$
while 
$$
d (\tau^* \psi_2^* \alpha_2)_{(x,0,0)}|_{\ker (\psi_1^* \alpha_1)_{(x,0,0)}} =
(id_\zeta \oplus id _{\zeta^*})^* \omega_{\zeta_x \oplus \zeta_x^*} 
\oplus \sigma ^* (\omega_2)_x 
= \omega_{\zeta_x \oplus \zeta_x^*} \oplus (e^h \omega_1)_x.
$$

We now apply the equivariant relative Darboux theorem to $N
\hookrightarrow \zeta^* \oplus E_1$, $\psi_1^* \alpha_1$ and $\tau^*
\psi_2^* \alpha_2$ to obtain a $G$-equivariant diffeomorphism $F: V_1
\to V_2$ with $F|_N = id_N$, $F^* \tau^* \psi_2^* \alpha_2 = e^f
(\psi_1^* \alpha_1)$ for some neighborhoods $V_1, V_2$ of the zero
section and some function $f$.  The theorem follows.
\end{proof}

\begin{lemma}\labell{lem.pre-iso.orb}
Let $\Psi _\alpha :M \to \fg^*$ be the $\alpha$-moment map for an
action of a torus $G$ on a contact manifold $(M, \alpha)$.  Suppose
for some point $x$ we have $\Psi_\alpha (x) \not =0$. Then the orbit
$G\cdot x$ is pre-isotropic in $(M, \xi = \ker \alpha)$.
\end{lemma}

\begin{proof}
Since the contact distribution $\xi$ is of codimension 1, in order to
show that the orbit $G\cdot x$ is transverse to $\xi$, it is enough to
prove that there is a vector $X\in \fg$ such that $X_M(x) \not \in
\xi_x$, i.e., such that 
$0 \not = \alpha _x (X_M (x)) = \langle \Psi_\alpha (x), X\rangle$. 
But the latter is exactly the condition that $\Psi_\alpha (x) \not = 0$. 

Next note that  the fiber of the characteristic distribution $\zeta$ at $x$
 is 
$$
\zeta_x  = T_x (G\cdot x) \cap \xi_x = \{ X_M (x) \mid 0= \alpha_x (X_M (x)) 
        = \langle \Psi _\alpha (x) , X\rangle \}.  
$$ 
Let $\fk = \{ X \in \fg \mid \langle \Psi_\alpha (x), X \rangle =
0\}$.  Since $\fg$ is abelian, $\fk$ is a subalgebra.  Consequently
$\zeta$ is an integrable distribution.  Since the leaves of the
foliation defined by $\zeta$ are tangent to the contact structure,
$\zeta$ is an isotropic subbundle of the symplectic vector bundle
$(\xi, \omega)$, where as usual $\omega = d\alpha |_\xi$.  
\end{proof}

\begin{definition} \labell{rmrk.sympl-slice}
Let $(M, \xi = \ker \alpha)$ be a contact manifold with an action of a
torus $G$ preserving the contact form $\alpha$ and let $\Psi_\alpha :M
\to \fg^*$ be the corresponding moment map.  Suppose a point $x\in M$ is such 
that $\Psi_\alpha (x) \not =0$.  Let $\zeta_x$ denotes the fiber at $x$ of 
the characteristic distribution of the pre-isotropic embedding $G\cdot x
\hookrightarrow (M, \xi)$ and let $\zeta_x^\omega$ denote its symplectic 
perpendicular in $(\xi, \omega = d\alpha |_\xi)$.  We define the {\bf
symplectic slice} at $x$ for the action of $G$ on $(M, \xi)$ to be the
conformal symplectic vector space $V = \zeta_x^\omega/ \zeta_x$ with
the conformal symplectic structure $[\omega_V]$ induced by
$[\omega]$. We refer to the symplectic representation of the isotropy
group $G_x$ on $(V, \omega_V)$ as the {\bf symplectic slice
representation}.

We define the {\bf characteristic subalgebra } of the embedding $G\cdot x
\hookrightarrow (M, \xi)$ to be $\fk := (\R \Psi_\alpha (x))^\circ$.  Note 
that $\zeta_x = \fk/\fg_x$ and that $\fk$ is co-oriented. 
\end{definition}

\begin{lemma}\labell{lem.unique-orb}
Let $(M_j, \xi_j = \ker \alpha_j)$, $j= 1,2$, be two contact manifolds
with actions of a torus $G$ preserving the contact forms $\alpha_1$,
$\alpha_2$.   Suppose $x_j \in M_j$, $j= 1,2$ are two points such that 
\begin{enumerate}
\item $0\not = \Psi_{\alpha_1} (x_1) = \lambda \Psi_{\alpha_2} (x_2)$ 
for some $\lambda >0$ (i.e., the characteristic subalgebras agree as
co-oriented subspaces of $\fg$);
\item the isotropy groups are equal : $G_{x_1} = G_{x_2}$;
\item the symplectic slice representations at $x_1$ and $x_2$ are isomorphic 
as symplectic representation up to a conformal factor.
\end{enumerate}
Then there exist $G$-invariant neighborhoods $U_j$ of $G\cdot x_j$ in
$M_j$, $j= 1, 2$, and a $G$-equivariant diffeomorphism $\varphi :U_1
\to U_2$ such that $\varphi^* \alpha_2 = e^f \alpha_1$ for some
function $f$.
\end{lemma}
\begin{proof}
The characteristic distributions and the conformal symplectic normal bundles of
the embeddings $\iota_j : G\cdot x_j \hookrightarrow (M_j, \xi_j)$,
$j=1,2$, are, respectively, 
$$
\zeta_j =  G \times _{G_{x_j}} (\fk/\fg_{x_j}) \quad \text{ and } \quad 
E_j =G \times _{G_{x_j}} V_j, 
$$ 
where $\fk = (\R \Psi_{\alpha_1} (x_1))^\circ = (\R \Psi _{\alpha_2}
(x_2))^\circ$ and $G_{x_j} \to \text{Sp} (V_j, \omega_{V_j})$ are the 
symplectic slice representations.
The lemma follows from the uniqueness of pre-isotropic embeddings 
(Theorem~\ref{thm-unique-pre}).
\end{proof}

\begin{lemma}\labell{lemma.pre-iso-model}
Let $(M, \xi = \ker \alpha)$ be a contact manifold with action of a
torus $G$ preserving the contact form $\alpha$.  Suppose $x\in M$ is
such that $\Psi_\alpha (x) \not = 0$.  Let $\fk = (\R\Psi_\alpha
(x))^\circ$ be the characteristic subalgebra and $G_x \to \text{Sp}
(V, \omega_V)$ the symplectic slice representation.  Choose splittings
$$
\aligned
\fg_x^\circ & = (\fk/\fg_x)^* \oplus \R \Psi_\alpha (x)  \cr
\fg^* &= \fg_x^\circ \oplus \fg_x^*
\endaligned
$$
and thereby a splitting
$$
\fg^* = (\fk/\fg_x)^* \oplus\R \Psi_\alpha (x)\oplus \fg_x^*.
$$ 
Let $i:\fg_x^* \hookrightarrow \fg^*$, $j:
(\fk/\fg_x)^*\hookrightarrow \fg^*$ be the corresponding embeddings.

There exists a $G$-invariant neighborhood $U$ of the zero section
$G\cdot [1, 0, 0]$ in 
$$
N = G\times _{G_x} ((\fk/\fg_x)^* \oplus V)
$$
and an open $G$-equivariant embedding $\varphi : U \hookrightarrow M$
with $\varphi([1, 0, 0]) = x$ and a $G$-invariant 1-form $\alpha_N$ on $N$ 
such that 
\begin{enumerate}
\item $\varphi ^*\alpha = e^f \alpha_N $ for some function 
$f\in C^\infty (U)$ and 

\item the $\alpha_N$-moment map $\Psi_{\alpha_N}$ is given by 
$$ 
\Psi_{\alpha_N}([a,\eta, v]) =
\Psi_\alpha (x) + j (\eta) + i (\Phi _V (v))
$$
where $\Phi_V :V \to \fg^*$ is the homogeneous moment map for the
slice representation.
\end{enumerate}
Consequently, 
$$
\Psi_\alpha \circ \varphi ([a,\eta, v]) = 
	(e^f \Psi_{\alpha_N})([a,\eta, v]) =
	e^{f([a, \eta, v])} (\Psi_\alpha (x) + j (\eta) + i (\Phi _V (v))),
$$ 
for some $G$-invariant function $f $ on $N$.
\end{lemma}

\begin{proof}
By Lemma~\ref{lem.unique-orb} it is enough to construct on $N= G\times
_{G_x} ((\fk/\fg_x)^* \times V)$ a $G$-invariant contact form
$\alpha_N$ so that the embedding $\iota : G/G_x \hookrightarrow (N,
\ker \alpha _N)$, $\iota (aG_x) = [a, 0, 0]$ is pre-isotropic, 
the symplectic slice representation at $[1,0,0]$ is $G_x \to \text{Sp}
(V, \omega_V)$ and $\Psi_{\alpha_N} ([1,0,0]) = \Psi_\alpha (x)$.

We construct $(N, \alpha_N)$ as a contact quotient 
(see Lemma~\ref{lemma-contact-quotient}).
Since $G$ is abelian, both right and left trivializations identify the
cotangent bundle $T^* G$ with $G\times \fg^*$.  Consider the hypersurface
$$
\Sigma = 
G\times \left ( \Psi_\alpha (x) + j ((\fk/\fg_x)^*) + i (\fg_x^*)\right) 
$$
in $G\times \fg^* = T^*G$.
Since $\Psi_\alpha (x) \not = 0$, $\Sigma $ is a hypersurface of
contact type  (the expanding vector field $X$ is
generated by dilations $(0, \infty)\times G\times \fg^* \ni (t, g,
\nu)\mapsto (g, t\nu)$).  Consider  the action of $G$ on 
$G\times \fg^*$ given by
 $g\cdot (a, \nu) = (ga, \nu)$ and the action of $G_x$ on $G\times
\fg^*$ given by $b\cdot (a, \nu) = (ab\inv, \nu)$. Both actions preserve
$\Sigma$, $X$ and the tautological 1-form $\alpha_{T^*G}$.  The action
of $G_x$ on $V$ preserves the 1-form $\alpha_V = \iota (R)\omega_V$
where $R$ is the radial vector field on $V$. The diagonal action of
$G_x$ on $\Sigma \times V$ preserves the contact form $\left(\alpha
_{T^*G}|_\Sigma\right) \oplus \alpha_V$.  The corresponding moment map $\Phi
:\Sigma \times V \to \fg_x^*$ is given by
$$
\Phi ((a, \Psi_\alpha (x) + j (\eta) + i(\mu)), v) = -\mu + \Phi_V (v) 
$$
where $\Phi_V :V \to \fg_x^*$ is the $\alpha_V$-moment map.  Therefore the
 reduced space at zero for the action of $G_x$ is 
$$
N:=\Phi\inv (0) /G_x \simeq \{ (a, \eta, \mu, v)
 \in G \times (\fk/\fg_x)^* \times
\fg_x^* \times V \mid \mu = \Phi _V (v)\}/G_x =
G\times _{G_x} ( (\fk/\fg_x)^* \times V) $$ and $\alpha
_{T^*G}|_\Sigma \oplus \alpha _V$ descends to a $G$-invariant contact
form $\alpha_N$ on $N$.

Note that the moment map for the action of $G$ on $\Sigma\times V$
descends to the $\alpha_N$-moment map for the induced action of $G$ on
$N$. Hence it is given by the desired formula:
$$
\Psi_{\alpha_N} ([a, \eta, v]) = \Psi_\alpha (x) + j(\eta) + i (\Phi_V (v)).
$$
\end{proof}

We will need the following standard fact.
\begin{lemma}\labell{weighty-lemma} 
 A symplectic representation of a torus has well-defined weights.
\end{lemma}
\begin{proof}
Since the unitary group is the maximal compact subgroup of the
symplectic group, given a symplectic representation of a torus $\rho :
H\to \text{Sp} (V, \omega)$ there exists on $V$ an $H$-invariant
contact structure $J$ compatible with the symplectic form $\omega$.
We define the weights of $\rho$ to be the weights of the complex
representation $\rho : H \to \text{GL} (V, J)$.  Since any two
$H$-invariant complex structures on $V$ compatible with $\omega$ are
homotopic, the weights do not depend on the choice of $J$, i.e., they
are well-defined.
\end{proof}

The next lemma is taken from \cite{D}.
\begin{lemma}\labell{lemma.B}
If $\rho : H \to \text{Sp}(V, \omega)$ is a faithful symplectic
representation of a compact abelian group $H$ and if $2\dim H = \dim
V$ then $H$ is connected and the weights of $\rho$ form a basis of the
weight lattice $\Z_H^*$ of $H$.
\end{lemma}
\begin{proof}
Arguing as in the proof of Lemma~\ref{weighty-lemma} we may assume
that $\rho $ maps $H$ into the unitary group $U(V)$.  Since the
dimension of a maximal torus of $U(V)$ is $\frac{1}{2} \dim V$ and
since $\rho$ is faithful, $\rho $ maps the identity component of $H$
{\bf onto} a maximal torus of $U(V)$.  Since the centralizer of a
maximal torus of $U(V)$ is the torus itself, $H$ is connected.
Finally, the weights of a maximal torus $T$ of $U(V)$ form a basis of
the weight lattice of $T$.
\end{proof}
\noindent
Recall that $G$-orbits in contact toric $G$-manifolds are
pre-isotropic (Lemma~\ref{lem.pre-iso.orb} and
Lemma~\ref{lemma-nonzero}).

\begin{lemma}\labell{lemma.conn+slice}
Let $(M, \alpha, \Psi_\alpha :M \to \fg^*)$ be a contact toric
$G$-manifold.  For any point $x\in M$ the symplectic slice
representation $\rho : G_x \to \text{Sp}(V)$ is faithful and $\dim G_x
= \frac{1}{2} \dim V$.  Consequently the isotropy group $G_x$ is
connected.  Also the image of the moment map $\Phi_V (V)$ for the
slice representation $\rho$ has the following properties: the cone
$\Phi_V (V)$ has $d = \dim G_x$ edges; each edge is spanned by a
weight of $G_x$; these weights form a basis of the integral lattice of
$G_x$.  Hence the cone $\Phi_V (V)$ completely determines the slice
representation $\rho$.
\end{lemma}

\begin{proof}
Let $\fk = (\R \Psi_\alpha (x))^\circ$ denote the characteristic
subalgebra of the pre-isotropic embedding $G\cdot x \hookrightarrow
(M, \xi =\ker \alpha)$.  By Lemma~\ref{lemma.pre-iso-model} a
neighborhood of $G\cdot x$ in $M$ is $G$-equivariantly diffeomorphic
to a neighborhood of the zero section in $N = G\times _{G_x}
((\fk/\fg_x)^* \times V)$.  Since $G$ is abelian, the action of $G_x$
on $(\fk/\fg_x)^*$ is trivial.  Since by assumption the action of $G$
on $M$ is effective, the slice representation of $G_x$ on $V$ has to
be faithful.

By definition of $V$, the dimension of $V$ is the dimension of the
contact distribution minus twice the dimension of the characteristic
distribution, i.e., $\dim V = (\dim M -1) - 2 \dim (\fk/\fg_x)$.  Now
$\dim M - 1 = 2 \dim G -2$ and $\dim \fk = \dim G -1$.  Therefore,
$\dim V = 2\dim G -2 - 2 ((\dim G -1) - \dim G_x) = 2 \dim G_x$.  By
Lemma~\ref{lemma.B} $G_x$ is connected and the weights $\nu_1, \ldots
\nu_d$ ($d = \dim G_x = \frac{1}{2} \dim V$)  of the slice representation 
$\rho$ form a basis of the weight lattice of $G_x$. On the other hand,
$\Phi_V (V) = \{ \sum _{i=1}^d a_i \nu_i \mid a_i\geq 0\}$.  The rest
of the lemma  follows.
\end{proof}

As a corollary of the first part of the proof and
Lemma~\ref{lemma.pre-iso-model} we get
\begin{theorem}\label{BlockF}
Let $(M, \alpha, \Psi_\alpha: M \to \fg^*)$ be a c.c.c.t. $G$-manifold
normalized so that 
$\Psi_\alpha (M)\subset S(\fg^*) = 
\{ \eta \in \fg^* \mid ||\eta || = 1\}$.  Let $x\in M$ be
a point, $G_x$ be its isotropy group (which is connected). Let $\rho:
G_x \to \text{Sp}(V, \omega_V)$ denote the symplectic slice
representation, $\Phi: V \to \fg^*$ denote the corresponding moment
map, and let $\fk = (\R \Psi_\alpha (x))^\circ$ be the characteristic
subalgebra. Choose the embeddings $i:\fg_x^* \to \fg^*$, $j:
(\fk/\fg_x)^* \to \fg^*$ as in Lemma~\ref{lemma.pre-iso-model}.

There exists an open embedding $\varphi$ from a neighborhood of the orbit 
$G/G_x \times \{0\}\times \{0\}$ in $G/G_x \times (\fk/\fg_x)^* \times V$ into
$M$ such that
\begin{equation}
\left(\Psi_\alpha \circ \varphi\right) (aG_x, \eta, v) 
= \frac{\Psi_\alpha (x) + j (\eta) + i (\Phi (v))}
{||\Psi_\alpha (x) + j (\eta) + i (\Phi (v))||}.
\end{equation}
\end{theorem}
\begin{proof}
Since the isotropy group $G_x$ is connected, the sequence $1\to G_x
\to G \to G/G_x \to 1$ splits.  Hence $G\times _{G_x} ((\fk/\fg_x)^*
\times V) = G/G_x \times (\fk/\fg_x)^* \times V$.
\end{proof}

\begin{remark}\labell{rmrk.slice-moment}
Recall that for the standard representation of the $n$-torus $\bbT ^n$
on $\C ^n$ preserving the standard symplectic form $\omega = \sqrt{-1}
\sum dz_j \wedge d \bar{z}_j$, the corresponding moment map $\Phi: \C^n
\to (\R^n)^*$ is given by $\Phi (z) = (|z_1|^2, \ldots , |z_n|^2)$.  Hence 
the fibers of $\Phi$ are $\bbT ^n$-orbits.  Consequently if $\rho : H
\to \text{Sp}(V, \omega_V)$ is a faithful representation of a torus
$H$ with $\dim V = 2 \dim H$, then the fibers of the corresponding
moment map $\Phi_V : V \to \fh^*$ are $H$-orbits.
\end{remark}

\begin{lemma}\labell{BlockC1}
Let $(M, \alpha, \Psi_\alpha: M \to \fg^*)$ be a compact connected
contact toric $G$-manifold.
Then
\begin{enumerate}
\item The connected components of the fibers of $\Psi_\alpha$ are $G$-orbits.
\item For any point $x\in M$ and any sufficiently small $G$-invariant 
neighborhood $U$ of $x$ in $M$ the pair $(\R^+ \Psi_\alpha (x),
C(\Psi_\alpha |_U))$ determines the contact toric manifold $(U, \alpha
|_U, \Psi_\alpha |_U = \Psi _{\alpha |_U})$. (Recall that 
$C(\Psi_\alpha |_U)) = 
\{ t \Psi_\alpha (x) \mid t\in [0, \infty), \, x \in U\}$ 
is the moment cone of $\Psi_\alpha$).
\end{enumerate} 
\end{lemma}
\begin{proof}
Fix a point $x\in M$.  Let $\rho: G_x \to \text{Sp}(V, \omega_V)$
denotes the corresponding symplectic slice representation and $\fk =
(\R \Psi_\alpha (x))^\circ$ the characteristic subalgebra.  Let $i:
\fg^*_x \to \fg^*$ and $j: (\fk/\fg_x)^* \to \fg^*$ be the embeddings 
as in Lemma~\ref{lemma.pre-iso-model}.  

By Lemma~\ref{lemma.conn+slice}, the isotropy group $G_x$ is
connected.  By Lemma~\ref{lemma.pre-iso-model} there exists a
$G$-invariant neighborhood $U$ of $G\cdot x$ in $M$, a $G$-invariant
neighborhood $U_0$ of $G\cdot [1, 0,0]$ in $N = G\times _{G_x}
((\fk/\fg_x)^* \times V)$, a $G$-invariant contact form $\alpha_N$ on
$N$ and a $G$-equivariant diffeomorphism $\varphi : U_0 \to U$ such
that $\varphi^* \alpha = e^f\alpha_N$ for some $G$-invariant function
$f$.  Consequently the $\alpha$- and $\alpha_N$-moment maps are
related by $\Psi_\alpha \circ \varphi ([a, \eta, v]) 
= e^{f ([a, \eta, v])}\Psi_{\alpha_N} ([a,\eta, v])$.
Recall that the $\alpha_N$-moment map $\Psi_{\alpha_N}$ is given
by $\Psi_{\alpha_N} ([a, \eta, v]) = (\Psi_\alpha (x) + j (\eta) + i
(\Phi_V (v)))$, where $\Phi_V : V
\to \fg_x^*$ is the moment map for the slice representation.  

As observed previously the connectedness of $G_x$ implies that $N $ is
diffeomorphic to $G/G_x \times (\fk/\fg_x)^*\times V$.  Under this
identification $\Psi_{\alpha_N} (aG_x, \eta, v) =
\Psi_\alpha (x) + j (\eta) + i (\Phi_V(v))$.  Hence a fiber of
$\Psi_{\alpha_N}$ is of the form $G/G_x \times \{\eta\} \times \Phi_V \inv
(\mu)$.  Since the fibers of $\Phi_V$ are $G_x$-orbits (see
Remark~\ref{rmrk.slice-moment}), the fibers of $\Psi_{\alpha_ N}$ are
$G$-orbits.  It follows that for any point $\eta \in \fg^*$ the set
$\Psi_\alpha \inv (\eta) \cap U$ is a $G$-orbit.

We next argue that the pair $( \R^+ \Psi_\alpha (x), C(\Psi_\alpha
|_U))$ determines $G_x$ and the symplectic slice representation $\rho$
(it obviously determines the characteristic subalgebra).  It is no
loss of generality to assume that $(M, \alpha, \Psi_\alpha) = (N,
\alpha_N, \Psi_{\alpha_N})$ (note that the contact form $\alpha_N$ is
not normalized but this won't matter) and that $U$ is a neighborhood
of $G/G_x \times \{0\} \times \{0\}$ in $N = G/G_x \times
(\fk/\fg_x)^* \times V$.  Since $G_x$ is connected it is determined by
its Lie algebra $\fg_x$ or, equivalently, by its annihilator
$\fg_x^\circ$.  Let $C$ be the moment cone of $(U, \alpha_N |_U,
\Psi_{\alpha_N}|_U)$. We may assume that $U$ is of the form $G/G_x
\times D_1 \times D_2$ where $D_1$ is a neighborhood of $0$ in
$(\fk/\fg_x)^*$ and $D_2$ is a $G_x$-invariant neighborhood of 0 in
$V$.  Then $C = \{0\} \cup \R^+ (\Psi_\alpha (x) + j (D_2) + i (\Phi_V
(D_2)))$.  Note that $\R^+ (\Psi_\alpha (x) + j (D_2))$ is an open
subset of $\fg_x^\circ$.  Hence $$
\fg_x^\circ = 
\{ w\in \fg^* \mid \text{there is } \epsilon > 0 \text{ such that } 
\Psi_\alpha(x) + t w \in C \text{ for all } t\in (-\epsilon, \epsilon )\}.
$$ 

Once we determined $\fg_x$ we have the natural inclusion $\iota:
\fg_x \to \fg$ and the dual projection $\iota^T : \fg^* \to \fg_x^*$.  
Note that $\iota^T \circ i = id_{\fg_x^*}$.  Therefore $\iota^T (C) =
\Phi_V (V)$.  By Lemma~\ref{lemma.conn+slice}, the cone $\Phi_V (V)$
completely determines the representation $\rho : G_x \to \text{Sp} (V,
\omega_V)$.
\end{proof}

\section{Properties of contact moment maps}

We need two properties of moment maps for action of tori on contact manifolds. 

\begin{definition}
Let $\Psi_\alpha : M \to \fg^*$ be the moment map for an action of a
torus $G$ on a manifold $M$ preserving a contact form $\alpha$.  The
corresponding {\bf orbital moment map} is the map $\bar{\Psi}_\alpha :
M/G \to \fg^*$ induced by $\Psi_\alpha$.
\end{definition}

The first property is the convexity of the image and the connectedness
of the fibers which is due to Banyaga and Molino in the toric case
\cite{BM1,BM2} (by Lemma~\ref{lemma-nonzero} we know that  
Theorem~\ref{main-convexity-theorem} below applies to contact toric
manifolds).\footnote{I don't
understand the proof of convexity and connectedness in \cite{BM2}.  In
particular it is not  clear to me how the hypothesis that the
dimension of the group is at least 3 is being used. }

\begin{theorem}\labell{main-convexity-theorem}
Let $(M, \xi)$ be a co-oriented contact manifold with an effective
action of a torus $G$ preserving the contact structure and its
co-orientation. Let $\xi^\circ_+$ be a component of the annihilator of
$\xi$ in $T^*M$ minus the zero section: $\xi^\circ\smallsetminus 0
=\xi^\circ_+ \sqcup (-\xi^\circ_+ )$.  Assume that $M$ is compact and
connected and that the dimension of $G$ is bigger than 2.  If 0 is not
in the image of the contact moment map $\Psi: \xi^\circ_+ \to \fg^*$
then the fibers of $\Psi$ are connected and the moment cone $C(\Psi) =
\Psi (\xi ^\circ_+)
\cup \{0\}$ is a convex rational polyhedral cone.
\end{theorem}
\begin{proof}
See \cite{L-convex}.
\end{proof}

\begin{lemma}\labell{lemma-*} 
Let $(M, \alpha, \Psi_\alpha: M \to \fg^*)$ be a c.c.c.t. $G$-manifold
normalized so that $\Psi_\alpha (M)\subset S(\fg^*)$.  Suppose $\dim M
> 3$.  Then the fibers of the moment map $\Psi_\alpha$ are
$G$-orbits. Consequently the orbital moment map $\bar{\Psi}_\alpha
:M/G \to S(\fg^*)$ is a (topological) embedding.
\end{lemma}

\begin{proof}
By Theorem~\ref{main-convexity-theorem}, the fibers of $\Psi_\alpha$
are connected.  By Lemma~\ref{BlockC1} the connected components of the
fibers of $\Psi_\alpha$ are $G$-orbits.  Therefore the fibers of
$\Psi_\alpha$ are $G$-orbits and the orbital moment map is injective.
Since $M/G$ is compact, the orbital moment map $\bar{\Psi}_\alpha$ is
an embedding.
\end{proof}

Lemma~\ref{lemma-*} does not hold for 3-dimensional contact toric
manifolds.  For example consider $(M, \alpha) = (S^1 \times \bbT^2,
\cos nt\, d\theta_1 + \sin nt \, d\theta _2)$ with the corresponding
moment map $\Psi_\alpha (t, \theta_1, \theta_2) = (\cos nt, \sin nt)$.
The orbital moment map is $\bar{\Psi}_\alpha : S^1 = M/\bbT^2 \to S(\R^2) = S^1
$ is given by $\bar{\Psi}_\alpha (t) = (\cos nt, \sin nt)$, which is
not an embedding for $n>1$.  However, $\bar{\Psi}_\alpha (t)$ is an
embedding {\bf locally}.  This is true in general.

\begin{lemma}\labell{lemma-local-embedd}
Let $(M, \alpha, \Psi_\alpha: M \to \fg^*)$ be a c.c.c.t. $G$-manifold
normalized so that $\Psi_\alpha (M)\subset S(\fg^*)$. For any $x\in
M/G$ there is a neighborhood $U$ of $x$ in $M/G$ such that the
restriction of the orbital moment map $\bar{\Psi}_\alpha$ to $U$ is an
embedding into $S(\fg^*)$.
\end{lemma}
\begin{proof}
This is an easy consequence of the local normal form theorem,
Theorem~\ref{BlockF}.
\end{proof}

\begin{lemma}\labell{preBlockC3}
Let $(M, \alpha, \Psi_\alpha: M \to \fg^*)$ be a c.c.c.t. $G$-manifold
normalized so that $\Psi_\alpha (M)\subset S(\fg^*)$.  Assume $\dim M
> 3$.  If the moment map $\Psi_\alpha :M \to S(\fg^*)$ is onto, then
the action of $G$ on $M$ is free, hence $\Psi_\alpha :M \to S(\fg^*)$
is a principal $G$-bundle.
\end{lemma}

\begin{proof}
Suppose the action of $G$ is not free.  Then for some point $x\in M$
the isotropy group $G_x$ is not trivial.  By Lemma~\ref{lemma-*} for
any $G$-invariant neighborhood $U$ of the orbit $G\cdot x$ there is an
open subset $\widetilde{W}$ of the sphere $S(\fg^*)$ such that 
$$
\Psi_\alpha (U) = \widetilde{W} \cap \Psi_\alpha (M).
$$ 
By the local normal form theorem, Theorem~\ref{BlockF}, we may
choose $U$ and $\widetilde{W}$ so that
$$
\Psi_\alpha (U) = \widetilde{W} \cap 
\left\{ \left.
\frac{\Psi_\alpha (x) + j (\eta) + i (\Phi (v))}
{||\Psi_\alpha (x) + j (\eta) + i (\Phi (v))||} \right | 
\eta \in (\fk/\fg_x)^*, \, v \in V
\right\}
$$ 
where as in Theorem~\ref{BlockF} $\rho: G_x \to \text{Sp}(V,
\omega_V)$ denotes the symplectic slice representation, $\Phi: V \to
\fg^*$ denote the corresponding moment map, $\fk = (\R \Psi_\alpha
(x))^\circ$ is the characteristic subalgebra, and $i:\fg_x^* \to
\fg^*$, $j: (\fk/\fg_x)^* \to \fg^*$ are the embeddings  as in
Lemma~\ref{lemma.pre-iso-model}.  Since $G_x$ is nontrivial, the
symplectic slice $V$ is not zero.  It follows from
Lemma~\ref{lemma.conn+slice} that $\Phi (V)$ is a proper cone in
$\fg_x^*$.  Therefore $\Psi_\alpha \cap \widetilde{W} \not =
\widetilde{W}$, i.e., $\Psi_\alpha$ is not onto.  Contradiction. 
Therefore the action of $G$ is free.

By Lemma~\ref{lemma-*}, the fibers of $\Psi_\alpha$ are $G$-orbits.
Therefore if the action of $G$ is free, then $\Psi_\alpha: M \to
S(\fg^*)$ is a principal $G$-bundle.
\end{proof}
The next lemma is a partial converse to Lemma~\ref{preBlockC3}.
\begin{lemma}\labell{BlockC4}
Suppose a Lie group $G$ acts on a manifold $M$ preserving a contact
form $\alpha$.  Let $\Psi_\alpha :M \to \fg^*$ denotes the
corresponding moment map.  Suppose the action of $G$ at a point $x$ is
free and the value $\mu$ of moment map at $x$ is non-zero.  Then
$\pi_\mu \circ d(\Psi_\alpha)_x :T_x M \to \fg^* / \R \mu$ is onto.
Here $\pi_\mu : \fg^* \to \fg^*/\R\mu$ is the obvious projection.
\end{lemma}

\begin{proof}
Since the action of $G$ at $x$ is free, for any $0\not = X \in \fg$
the induced vector field $X_M $ is nonzero at $x$: $X_M (x) \not = 0$.
Since the action of $G$ preserves $\alpha$, $0 = d \iota (X_M) \alpha
+ \iota (X_M) d\alpha$.  Therefore, for any $v\in T_x M$
$$
\langle d(\Psi_\alpha)_x (v), X\rangle = d \langle \Psi_\alpha, X\rangle _x (v)
 = d\alpha _x (v, X_M (x)).
$$
Now 
\begin{equation}
\begin{split}
(\fg^*/\R\mu)^* = \ker \mu 
	& = \{X\in \fg \mid \langle  \Psi_\alpha (x), X\rangle = 0\}\\
	&=  \{X\in \fg \mid \alpha _x (X_M (x)) =0\}\\
	&=  \{X \in \fg \mid X_M \in \ker \alpha_x\}. 
\end{split}
\end{equation}
Hence to prove that $\pi_\mu \circ d(\Psi_\alpha)_x : T_xM \to
\fg^*/\R\mu$ is onto, it's enough to show that for any $0\not =X\in \fg$ 
with $\alpha_x (X_M (x)) =0$ there is $v\in T_xM$ with $d\alpha _x (v,
X_M (x)) =\langle d(\Psi_\alpha)_x (v), X\rangle \not = 0$.  Since
$\alpha$ is contact, $d\alpha |_{\ker \alpha}$ is nondegenerate.
Therefore, for any $X_M (x) \not= 0$ with $\alpha_x (X_M (x)) = 0$
there is $v \in \ker \alpha_x $ so that $d\alpha_x (v, X_M (x)) \not =
0$.
\end{proof}

\begin{corollary}\labell{Cor1}
Let $(M, \alpha, \Psi_\alpha: M \to \fg^*)$ be a c.c.c.t. $G$-manifold
normalized so that $\Psi_\alpha (M)\subset S(\fg^*)$.  If the action
of $G$ is free then the moment map $\Psi_\alpha: M \to S(\fg^*)$ is a
submersion.  

If additionally $\dim M > 3$ then $\Psi_\alpha: M \to
S(\fg^*)$ is a principal $G$-bundle.
\end{corollary}
\begin{proof}
By Lemma~\ref{BlockC4}, the differential $d(\Psi_\alpha)_x : T_x M \to
T_{\Psi_\alpha (x)} S(\fg^*)$ is surjective for all $x \in M$.
Consequently the image $\Psi_\alpha (M)$ is open in the sphere.  On
the other hand the image is closed since $M$ is compact.  Thus the
image is the whole sphere and $\Psi_\alpha : M \to S(\fg^*)$ is a
submersion. Since $M$ is compact it follows that $\Psi_\alpha : M \to
S(\fg^*)$ is a fibration.  If $\dim M > 3$ then, by
Lemma~\ref{lemma-*}, $\Psi_\alpha : M \to S(\fg^*)$ is a principal
$G$-bundle.
\end{proof}



\begin{definition} Two contact toric $G$-manifolds 
$(M, \alpha, \Psi_\alpha)$ and $(M, \alpha', \Psi_{\alpha'})$ are {\bf
locally isomorphic } if 
\begin{enumerate}
\item there exists a homeomorphism $\bar{\varphi}: M/G \to M'/G$ and
\item for any point $x\in M/G$ there is a neighborhood $U \subset M/G$ 
containing it and an isomorphism of contact manifolds 
$\varphi_U : \pi \inv (U) \to
(\pi')\inv (\bar{\varphi} (U))$ such that $\pi' \circ \varphi_U =
\bar{\varphi}\circ \pi$ where $\pi: M \to M/G$ and $\pi': M \to M/G$
are the orbit maps (hence, by Remark~\ref{rmrk2.15}, $\varphi_U^*
\alpha' = \alpha$).
\end{enumerate}
\end{definition}

\begin{lemma}\labell{BlockA}
Let $(M_1, \alpha_1, \Psi_{\alpha_1}: M_1\to\fg^*)$ and $(M_2,
\alpha_2, \Psi_{\alpha_2}: M_2\to\fg^*)$ be two c.c.c.t. $G$-manifolds 
normalized so that $\Psi_{\alpha_i} (M_i) \subset S(\fg^*)$, $i=1,2$.
Suppose there is a homeomorphism $\bar{\varphi}: M_1/G \to M_2/G$ so that
$\bar{\Psi}_{\alpha_2} \circ \bar{\varphi}= \bar{\Psi}_{\alpha _1}$, where
$ \bar{\Psi}_{\alpha _1}$, $\bar{\Psi}_{\alpha_2}$ are orbital moment maps.

Then $(M_1, \alpha_1, \Psi_{\alpha_1})$ and $(M_2,\alpha_2,
\Psi_{\alpha_2})$ are locally isomorphic.
\end{lemma}

\begin{proof}
Denote the orbit map $M_i \to M_i/G$ by $\pi_i$, $i=1,2$. We want to
show that for any point $x\in M_1$ there is a $G$-invariant
neighborhood $U_1 \subset M_1$ and a $G$-equivariant diffeomorphism
$\varphi_U: U_1 \to U_2 = \pi_2 \inv (\bar{\varphi} (\pi_1 (U_1)))$
such that $\varphi_U^* \alpha_2 =  \alpha_1$  and such that 
$$
\pi_2 \circ\varphi_U = \bar{\varphi} \circ (\pi_1 |_{U_1}).
$$

 Pick a point $x_2 \in \pi_2 \inv (\bar{\varphi} (x_1))$. Then,
since $\bar{\Psi}_{\alpha_2} \circ \bar{\varphi}= \bar{\Psi}_{\alpha
_1}$, we have $\Psi_{\alpha_2} (x_2) = \Psi_{\alpha_1} (x_1)$.  Also,
given a $G$-invariant neighborhood $U_1$ of $x_1$ in $M_1$, let $U_2 =
\pi_2 \inv (\bar{\varphi} (\pi_1 (U_1)))$.  We have $\Psi_{\alpha_2}
(U_2) = \bar{\Psi}_{\alpha_2} (\bar{\varphi} (\pi_1 (U_1)) =
\bar{\Psi} _{\alpha_1} (\pi_1 (U_1)) = \Psi_{\alpha_1} (U_1)$.
Therefore by Lemma~\ref{BlockC1} (2) if $U_1$ is sufficiently small
there exists a $G$-equivariant contactomorphism $\varphi_U :U_1 \to
U_2$.  Hence, as remarked earlier, $\Psi_{\alpha_2}
\circ \varphi_U =  \Psi_{\alpha_1}$ and $\varphi_U^*
\alpha_2 = \alpha_1$.  

It remains to show that the map $\bar{\varphi}_U$ induced by
$\varphi_U$ on $U_1/G$ is $\bar{\varphi}$.  If $U_1$ is sufficiently
small, then by Lemma~\ref{lemma-local-embedd} the maps
$\bar{\Psi}_{\alpha_i}: U_i/G \to S(\fg^*)$, $i=1,2$, are embeddings.
Since $\Psi_{\alpha_2}
\circ \varphi_U =  \Psi_{\alpha_1}$, $\bar{\Psi}_{\alpha_2}
\circ \bar{\varphi}_U =  \bar{\Psi}_{\alpha_1}$.  Therefore $\bar{\varphi}_U = 
\left(\bar{\Psi}_{\alpha_2}|_{\Psi_{\alpha_1} (U_1)}\right)\inv \circ 
\bar{\Psi}_{\alpha_1}|_{\pi_1(U_1)} = \bar{\varphi}|_{\pi_1(U_1)}$.
\end{proof}

\section{From local to global}

In this section we prove that compact connected contact toric
manifolds are classified by the elements of the first \v Cech
cohomology of their orbit space with coefficients in a certain sheaf.
The argument here is an adaptation of the argument in \cite{LT} (which
was due to Lerman, Tolman and Woodward), which, in tern, was an
adaptation of the argument in \cite{HS}.  I recently learned that
essentially the same idea was developed earlier by Boucetta and Molino
\cite{Bou-M}.

Let $(M, \xi = \ker \alpha )$ be a co-oriented contact manifold.  Recall that a
vector field $\Xi$ is {\bf contact} if its flow preserves the contact
distribution $\xi$, or, equivalently, if $L_\Xi \alpha = f\alpha$ for
some function $f\in C^\infty (M)$.  A choice of a contact form $\alpha$
with $\ker\alpha = \xi$ establishes a bijection between contact vector
fields and functions: given a contact vector field $\Xi$ the
corresponding function is $\alpha (\Xi)$.  Conversely, given a
function $f\in C^\infty (M)$ the corresponding contact vector field
$\Xi _f$ is defined by
\begin{equation}
\Xi _f = f Y_\alpha - (d\alpha |_\xi) \inv (df|_\xi),
\end{equation}
where $Y_\alpha$ is the Reeb vector field of $\alpha$, that is, the
unique vector field such that $\alpha (Y_\alpha) = 1$, $\iota
(Y_\alpha) d\alpha = 0$.  Note that since $d\alpha |_\xi $ is
nondegenerate, $(d\alpha |_\xi) \inv (df|_\xi)$ is a well-defined
vector field tangent to $\xi$. Also, if $\Xi $ is contact then 
$$
L_\Xi \alpha = Y_\alpha (\alpha (\Xi)) \alpha.
$$

If a Lie group $G$ acts on the manifold $M$ and if $\alpha $ is
$G$-invariant, then it induces a bijection between $G$-invariant
contact vector fields and $G$-invariant functions.

\begin{lemma}\labell{lemma?}
Suppose $(M, \alpha, \Psi_\alpha :M\to \fg^*)$ is a contact toric
$G$-manifold.  For any $G$-invariant function $f$ the flow
$\varphi^f_t$ of the corresponding contact vector field $\Xi _f $
preserves the contact form $\alpha$ and induces the identity map on
the orbit space $M/G$.  In particular $\Xi _f$ is tangent to
$G$-orbits.
\end{lemma}

\begin{proof}
We first consider the special case of $f=1$.  Then the corresponding
contact vector field is the Reeb vector field $Y_\alpha$.  By
definition of $Y_\alpha$ the Lie derivative $L_{Y_\alpha} \alpha$
satisfies $L_{Y_\alpha} \alpha = d \iota (Y_\alpha )\alpha + \iota
(Y_\alpha) d \alpha = d (1) + 0 = 0 $.  Hence for any $X \in \fg$ 
$$
L_{Y_\alpha} \langle \Psi _\alpha, X\rangle = L_{Y_\alpha} (\iota (X_M
) \alpha) = \iota (L_{Y_\alpha} X_M) \alpha + \iota (X_M
)(L_{Y_\alpha} \alpha) .  
$$ 

Since $Y_\alpha$ is unique, it is $G$-invariant.  Hence $[Y_\alpha,
X_M] = 0$ for any $X\in \fg$.  And by the previous computation
$L_{Y_\alpha} \alpha =0$.  Therefore 
\begin{equation}\labell{eq5.1}
L_{Y_\alpha} \langle \Psi
_\alpha, X\rangle = L_{Y_\alpha} (\iota (X_M) \alpha) = \iota
([Y_\alpha, X_M]) \alpha + \iota (X_M) L_{Y_\alpha} \alpha = 0+ 0.
\end{equation}
Thus the Reeb vector field is tangent to the fibers of the moment map
$\Psi_\alpha$. Since $(M, \alpha, \Psi_\alpha :M\to \fg^*)$ is toric,
the connected components of the fibers of $\Psi_\alpha$ are $G$-orbits
(see Lemma~\ref{BlockC1}).  Therefore the Reeb vector field is tangent to
$G$-orbits.  Hence for any $G$-invariant function $f$, $Y_\alpha (f) =
0$ and consequently $L_{\Xi _f}\alpha = Y_\alpha (f) \alpha =0$.  That
is, the flow $\varphi^f_t$ of $\Xi _f$ preserves $\alpha$.

 For any $X \in \fg$
$$
L_{\Xi_f} \langle \Psi _\alpha, X\rangle = L_{\Xi_f} (\iota (X_M
) \alpha) = \iota ([\Xi_f, X_M]) \alpha + \iota (X_M
)(L_{\Xi_f} \alpha)= 0+ 0.
$$
Since connected components of the fibers of $\Psi_\alpha$ are
$G$-orbits, this proves that the contact vector field $\Xi_f$ is
tangent to $G$-orbits for any invariant function $f$.  Hence the flow
of $\Xi_f$ induces the identity map on the orbit space $M/G$.
\end{proof}

\begin{proposition} \labell{prop1*}
For a fixed torus $G$, the isomorphism classes of contact toric
$G$-manifold locally isomorphic to a given contact toric $G$-manifold
$(M, \alpha, \Psi_\alpha)$ are in one-to-one correspondence with the
elements of the first \v Cech cohomology group $H^1 (M/G, \cS)$ where
$\cS$ is the sheaf of groups on the orbit space $M/G$ defined by
$$
\cS (U) =\Iso (\pi\inv (U)),
$$ 
the group of isomorphisms of the contact toric manifold $(\pi\inv
(U), \alpha |_{\pi\inv (U)}, \Psi_\alpha |_{\pi\inv (U)})$ (cf.\
Definition~\ref{def2.15}). Here again $\pi : M \to M/G$ denotes the
orbit map.
\end{proposition}

\begin{proof}
The argument is standard (compare \cite{HS}, Proposition~4.2 or \cite{Bou-M}).
Suppose $(M', \alpha', \Psi_{\alpha'})$ is a contact toric
$G$-manifold locally isomorphic to $(M, \alpha, \Psi_\alpha)$.  Fix a
homeomorphism $\bar{\varphi}: \bar{M} = M/G \to \bar{M}' = M'/G$.
Choose an open cover $\{V_i\}$ of $\bar{M}$ such that for each $i$
there is a $G$-equivariant contact diffeomorphism $\sigma_i : \pi\inv
(V_i) \to (\pi')\inv (\bar{\varphi} (V_i))$ inducing $\bar{\varphi}$
on $V_i$ (here $\pi' : M' \to M'/G$ is the orbit map).  Let 
$$
f_{ij} = \sigma_i \circ \sigma_j|_{V_i \cap V_j}; 
$$
it is a \v Cech 1-cocycle whose cohomology class in $H^1 (\bar{M},
\cS)$ is independent of the choices made to define it.

Conversely, given an element of $H^1 (\bar{M}, \cS)$ we can represent
it by a \v Cech cocycle $\{f_{ij} :\pi\inv (V_i \cap V_j) \to \pi\inv
(V_i \cap V_j)\}$. We construct the corresponding contact toric
$G$-manifold by taking the disjoint union of the manifolds $(\pi\inv
(V_i), \alpha |_{\pi\inv (V_i)}, \Psi_\alpha |_{\pi\inv (V_i)})$ and
gluing $\pi\inv (V_i)$ to $\pi\inv(V_j)$ along $\pi \inv(V_i\cap V_j)$
using $f_{ij}$.  The cocycle condition guarantees that the gluing is
consistent.
\end{proof}

\begin{proposition}\labell{prop4.3}
Let $(M, \alpha, \Psi_\alpha )$ be a contact toric $G$-manifold. Let
$\pi :M \to M/G$ denote the orbit map, and let $\Z_G := \ker \{ \exp :
\fg \to G\}$ denote the integral lattice of the torus $G$.  There
exists a short exact sequence of sheaves of groups
\begin{equation}\labell{eq.ses}
0 \to \underline{\Z}_G \stackrel{j}{\to} \cC \stackrel{\Lambda}{\to} \cS \to 1,
\end{equation}
where for a sufficiently small  open subset $U$ of the orbit space $M/G$
\begin{enumerate}
\item $\Z_G(U) := C^\infty (\pi\inv (U), \Z_G)^G$;
\item $\cC (U) = C^\infty (\pi \inv (U))^G$, the sheaf of ``smooth'' 
functions on $M/G$;
\item $\cS (U) := \Iso (\pi\inv (U))$ is the sheaf  defined in 
Proposition~\ref{prop1*}.
\end{enumerate}
Hence $\cS$ is a sheaf of abelian groups and the cohomology groups
$H^i (M/G, \cS)$ are defined for all indices $i \geq 0$.
\end{proposition}

\begin{proof}
Let $U$ be an open subset of $\bar{M} = M/G$ and let $f\in\cC (U)$.
By Lemma~\ref{lemma?}, the time $t$-flow $\varphi_t^f : \pi\inv (U)
\to \pi\inv (U)$ induces the identity map on $U$ and preserves the
contact form $\alpha$.  Hence for any $t$, $\varphi_t^f$ is in $\cS
(U)$.  We define the map $\Lambda : \cC (U) \to \cS (U)$ by 
$$
\Lambda (f) = \varphi ^f _1.
$$
 
We next argue that $\Lambda$ is onto. Suppose $\varphi \in \cS (U)$.
By Theorem~3.1 of \cite{HS},
there exists a smooth $G$-invariant map $\sigma : \pi\inv (U) \to G$ such that
$$
\varphi (x) = \sigma (x)\cdot x
$$ 
for all $x\in \pi\inv (U)$.  Moreover, if $U$ is contractible, then
$\sigma (x)= \exp (X(x))$ for some smooth $G$-invariant map 
$X :\pi\inv (U) \to \fg$.  It's not hard to check that 
$x \mapsto \exp (X(x))\cdot x$ is the time-1 flow of the vector field
$\tilde{X}(x) := (X(x)_M) (x)$.  Note that $\tilde{X}$ is a
$G$-invariant vector field: for any $g\in G$ and $x\in \pi\inv (U)$
$$
\tilde{X} (g\cdot x) = (X(g\cdot x))_M (g\cdot x) = 
\left. \frac{d}{ds}\right|_{s = 0} \exp (s X(g\cdot x)) \cdot (g\cdot x) =
\left. \frac{d}{ds}\right|_{s = 0} g \cdot (\exp (s X( x)) \cdot  x) = 
dg_M (\tilde{X }(x)),
$$ 
where the second equality holds because $ {X}$ is
$G$-invariant and $G$ is abelian.

We next prove that the Lie derivative of $\alpha$ with respect to
$\tilde{X}$ is zero: $L_{\tilde{X}} \alpha = 0$.  To do this we recall a
few facts about basic forms \cite{Ko}.  Given an action of a compact
Lie group $G$ on a manifold $M$, a form $\beta$ is {\bf basic} if it
is $G$-invariant and if for any $X \in \fg$, the contraction
$\iota(X_M) \beta$ is zero (for zero forms we only require
invariance).  The set of basic forms is a subcomplex of the de Rham
complex of differential forms, i.e., if $\beta $ is basic then so is
$d\beta$.  Also, if $\varphi :M \to M$ is a $G$-equivariant map
inducing the identity on $M/G$ and $\beta$ is basic, then $\varphi^*
\beta = \beta$.  This is because it is a closed condition that holds
on the open dense subset of points of principal orbit type.

We claim that $L_{\tilde{X}} \alpha$ is basic for the action of the
torus $G$ on the manifold $\pi\inv (U)$. Note that $L_{\tilde{X}}
\alpha = d\iota (\tilde{X})\alpha + \iota (\tilde {X})d\alpha$.  Since
$\tilde{X}$ and $\alpha$ are $G$-invariant, $\alpha (\tilde{X})$ is
$G$-invariant hence basic. Therefore $d\iota (\tilde{X})\alpha$ is
basic.  Also, since $\tilde{X}$ and $\alpha$ are $G$-invariant, the
second term $ \iota (\tilde {X})d\alpha$ is $G$-invariant.  It remains
to show that for any $Y\in \fg$, $0 = \iota (Y_M)[ \iota (\tilde
{X})d\alpha]$.  Now for any $Y, Z \in \fg$ 
$$ 
d\alpha (Y_M, Z_M) = Y_M
(\alpha (Z_M)) - Z_M (\alpha (Y_M)) - \alpha ([Y_M, Z_M]) = 0 - 0 - 0,
$$ 
because $\alpha (Y_M)$, $\alpha (Z_M)$ are $G$-invariant functions
and because $[Y_M, Z_M] = - ([Y,Z])_M = 0$ (since $G$ is abelian).
Therefore for any $x\in \pi \inv (U)$
$$
\iota (Y_M)[ \iota (\tilde{X})d\alpha] = d\alpha _x (Y_M (x), (X(x))_M (x)) =0.
$$

Next let $\tau_t$ denote the time-$t$ flow of the vector field
$\tilde{X}$. Clearly $\tau_t$ is $G$-equivariant and induces the identity map 
on $U$. Thus , since $(L_{\tilde{X}}
\alpha)$ is basic we have $(\tau_t) ^* (L_{\tilde{X}} \alpha)
=(L_{\tilde{X}} \alpha)$ for all $t$.  We also know that $\tau_1 =
\varphi $ and that $\varphi^* \alpha = \alpha$.   Therefore 
$$ 
0 = \varphi^* \alpha - \alpha = \tau_1 ^*\alpha - \tau_0 ^* \alpha 
= \int _0 ^1 \frac{d}{dt} (\tau_t^*
\alpha)\, dt  = \int _0^1 \tau_t^* (L_{\tilde{X}} \alpha) \, dt = 
\int _0^1 (L_{\tilde{X}} \alpha) \, dt = (L_{\tilde{X}} \alpha).
$$ We conclude that $\tilde{X}$ is a contact vector field.  We define
$f = \alpha (\tilde{X})$.  Then the contact vector field of $f$ is
$\tilde{X} $ and $\Lambda (f) = \varphi_1^f = \tau_1 = \varphi$.  This
concludes the proof that $\Lambda : \cC \to \cS$ is onto.

We define the map $j: \cL \to \cC$ by $j(X) = \langle \Psi_\alpha , X
\rangle$.  Thus it remains to show that for any sufficiently small set
$U \subset M/G$ and any function $f\in \cC (U)$ if $\varphi^f_1 (x) =
x$ for all $x \in \pi\inv (U)$ then $$ f =\langle \Psi_\alpha , X
\rangle $$ for some $X\in \Z_G$. In fact it is enough to show that the
above equation holds on the open dense subset $\pi\inv (U_0)$ of
$\pi\inv (U)$ consisting of the points where the action of $G$ is
free.  Now the contact vector field $\Xi _f $ of $f$ is $G$-invariant
and is tangent to $G$-orbits (Lemma~\ref{lemma?}).  Therefore there
exists on $\pi\inv (U_0)$ a $G$-invariant smooth map $X: \pi\inv (U_0)
\to \fg$ such that 
$$
\Xi _f (x) = (X (x))_M (x).
$$
 Now the time--1 flow of $\tilde {X} (x) := (X (x))_M (x)$ is
$x\mapsto \exp (X(x))\cdot x$.  Thus if $\varphi^f_1 (x) = x$, then
$\exp (X(x))\cdot x = x$ for all $x\in \pi\inv (U_0)$.  Hence $\exp
(X(x)) = 1$ and so $X(x) \in \Z_G$ for all $x\in \pi\inv (U_0)$.
Therefore, since $\Z_G$ is discrete and $X$ is continuous and since we
may take $\pi\inv (U_0)$ to be connected, $X(x) = X$ for some fixed
vector $X\in \Z_G$.  It follows, since $\pi\inv (U_0)$ is dense in
$\pi \inv (U)$ that $\Xi_f (x) = X_M (x)$ for all $x\in \pi\inv (U)$.
Consequently $f = \alpha (\Xi _f ) = \alpha (X_M) = \langle
\Psi_\alpha, X\rangle$.
\end{proof}

\begin{corollary}  \labell{cor4.4}
Under the hypotheses of the proposition above,
$$ 
H^i (M/G, \cS) = H^{i+1} (M/G, \underline{\Z}_G)
$$ for all $i>0$.
\end{corollary}
\begin{proof}
The sheaf $\cC$ is a fine sheaf, so the long exact sequence in
cohomology induced by (\ref{eq.ses}) breaks up for $i>0$ into
isomorphisms $H^i (M/G, \cS) \simeq H^{i+1} (M/G, \underline{\Z}_G) $.
\end{proof}

\section{Proof of the main classification theorem}

\subsection{\protect Proof of Theorem~\ref{main-theorem} (1)}

Suppose $(M, \alpha, \Psi_\alpha)$ is a c.c.c.t. $G= \bbT^2$ manifold
and suppose the action of $G$ is free. Then the orbit space $M/G$ is a
1-dimensional compact connected manifold without boundary, hence is a
circle $S^1$.  Moreover, since any principal $\bbT^2$-bundle over
$S^1$ is trivial, $M$ is diffeomorphic to $S^1 \times \bbT^2 =
\bbT^3$.  It remains to show that the contact form $\alpha$ equals
$\alpha_n = \cos nt\, d\theta_1 + \sin nt \, d\theta_2$ for some
positive integer $n$.  

By Lemma~\ref{BlockC4} the orbital moment map $\bar{\Psi}_\alpha :S^1
= M/G \to S(\fg^*) = S^1$ is a submersion, hence is a covering map.
Therefore there is a homeomorphism $\bar{\varphi}: M/G \to S^1$ such
that $\bar{\Psi}_{\alpha_n} \circ \bar{\varphi} = \bar{\Psi}_\alpha$
where $\Psi_{\alpha_n}$ is the $\alpha_n$-moment map and $n$ is the
number of sheets in the cover $\bar{\Psi}_\alpha :S^1 \to S^1$.  By
Lemma~\ref{BlockA} $(M, \alpha, \Psi_\alpha)$ is locally isomorphic to
$(M, \alpha_n, \Psi_{\alpha_n})$.  It follows from
Proposition~\ref{prop1*}, Corollary~\ref{cor4.4} and the fact that
$H^2(S^1, \Z^2)= 0$ that $(M, \alpha, \Psi_\alpha)$ is isomorphic to
$(M, \alpha_n, \Psi_{\alpha_n})$.

\subsection{\protect Proof of Theorem~\ref{main-theorem} (2)}

\begin{lemma} \labell{lemma6.1}
Let $(M, \alpha, \Psi_\alpha)$ be a c.c.c.t. $G= \bbT^2$ manifold
normalized so that $\Psi_\alpha (M) \subset S(\fg^*)$ and suppose the
action of $G$ is {\bf not} free. Then
\begin{enumerate}
\item The orbit space $M/G$ is homeomorphic to the interval $[0,1]$.

\item The orbital moment map $\bar{\Psi}_\alpha :M/G \to S(\fg^*) = S^1$ 
lifts to an embedding $\tilde{\Psi}_\alpha :M/G \to \R$ so that $\Pi
\circ \tilde{\Psi}_\alpha = \bar{\Psi}_\alpha$ where $\Pi: \R \to S^1$
is the covering map $\Pi(t) = (\cos t, \sin t)$.

\item $\tilde{\Psi}_\alpha (M/G) = [t_1, t_2]$,
 and $\tan t_1, \tan t_2$ are rational numbers.

\item  If $(M', \alpha', \Psi_{\alpha'})$ is another c.c.c.t. $G$-manifold 
with $\tilde{\Psi}_{\alpha'} (M'/G) = \tilde{\Psi}_\alpha (M/G)$ then
$(M', \alpha', \Psi_{\alpha'})$ is isomorphic to $(M, \alpha,
\Psi_{\alpha})$.

\item Given $t_1, t_2\in \R$ with $0\leq t_1 <2\pi$, $t_1 < t_2$ and 
$\tan t_1, \tan t_2$ rational, there is a c.c.c.t. $G$-manifold 
$(M, \alpha, \Psi_\alpha)$ with $\tilde{\Psi}_\alpha (M/G) = [t_1, t_2]$.
\end{enumerate}
\end{lemma}

\begin{proof}
Since $(M, \alpha, \Psi_\alpha)$ is contact toric, $\Psi_\alpha (x)
\not =0$ for any $x\in M$ (Lemma~\ref{lemma-nonzero}).  
Therefore the action of $G$ on $M$ has no
fixed points. By Lemma~\ref{lemma.conn+slice} all the isotropy groups
are connected.  Therefore they are either trivial or circles.  Suppose
the isotropy group $G_x$ at $x$ is a circle.  Then its Lie algebra $\fg_x$
equals the characteristic subalgebra $\fk$.  Hence $\Psi_\alpha (x)$
is a multiple of a weight $\mu \in \Z_G^*$.

Also, the dimension of the symplectic slice $V$ at $x$ is 2.  Hence
the symplectic slice representation is isomorphic to the standard
action of $S^1$ on $\C$: $\lambda \cdot z = \lambda z$.  Consequently,
by the local normal form theorem, Theorem~\ref{BlockF}, a neighborhood
of $x$ in $M$ is diffeomorphic to $S^1 \times \C$, and a neighborhood
of $G\cdot x$ in $M/G$ is homeomorphic to $\C/S^1 = [0, \infty)$.  We
conclude that $M/G$ is a 1-dimensional $C^0$ manifold with boundary.  Since
$M$ is compact and connected and since $G= \bbT^2$ it follows that
\begin{enumerate}
\item $M/G$ is homeomorphic to $[0, 1]$;

\item there are exactly two orbits $G\cdot x_1$, $G\cdot x_2$ which are 
diffeomorphic to $S^1$;

\item For $i=1,2$, $\Psi_\alpha (x_i) = \frac{\mu_i}{||\mu_i||}$ 
where $\mu_i = d\chi_i$ and the character $\chi_i$ is the map $G\to G/G_{x_i}$.
\end{enumerate}

Since $G=\bbT^2$ we may identify $\fg^*$ with $\R ^2$ and the weight
lattice $Z_G^*$ with $\Z^2$. Then for any weight $\mu$ of $G$,
$\frac{\mu}{||\mu||} = (\cos t, \sin t)$ for some $t\in \R$ with $\tan
t $ rational.

Since $M/G$ is homeomorphic to $[0, 1]$, the orbital moment map
$\bar{\Psi}_\alpha : M/G \to S(\fg^*) = S^1$ lifts to a map
$\tilde{\Psi}_\alpha :M/G \to \R$ such that $\Pi \circ
\tilde{\Psi}_\alpha = \bar{\Psi}_\alpha$.  Note that 
$\bar{\Psi}_\alpha ( M/G)$ is an interval with end points being the
images of the exceptional orbits $G\cdot x_1$, $G\cdot x_2$.  Hence we
may assume that $\bar{\Psi}_\alpha ( M/G) = [t_1, t_2]$, $0\leq t_1 <
2\pi$, and that $\tan t_1$, $\tan t_2$ are rational numbers.

Moreover since $\bar{\Psi}_\alpha$ is locally an embedding
(Lemma~\ref{lemma-local-embedd}), $\tilde{\Psi}_\alpha$ is an
embedding.  Thus $\tilde{\Psi}_\alpha :M/G \to [t_1, t_2]$ is a
homeomorphism.

Now suppose $(M', \alpha', \Psi_{\alpha'})$ is another
c.c.c.t. $G$-manifold with $\tilde{\Psi}_{\alpha'} (M'/G) = [t_1,
t_2]$.  Then $\varphi := (\tilde{\Psi}_{\alpha'})\inv \circ
\tilde{\Psi}_\alpha : M/G \to M'/G$ is a homeomorphism with
$\bar{\Psi}_{\alpha'} \circ \varphi = \bar{\Psi}_\alpha$.  By
Lemma~\ref{BlockA} $(M, \alpha, \Psi_\alpha)$ and $(M', \alpha',
\Psi_{\alpha'})$ are locally isomorphic.  Since $M/G$ is contractible, 
$ H^2 (M/G, \Z^2) =0$.  By Corollary~\ref{cor4.4} $H^1 (M/G, \cS) =
H^2 (M/G, \Z^2)$, where $\cS$ is the sheaf in
Proposition~\ref{prop1*}.  Hence, by Proposition~\ref{prop1*}, $(M,
\alpha, \Psi_\alpha)$ and $(M', \alpha',
\Psi_{\alpha'})$ are  isomorphic.

To prove part (5) we use an equivariant version of Proposition~2.15 in
\cite{L-IJM}:

\begin{proposition}\labell{cont-prop2}
Suppose $(\tilde{M}, \alpha)$ is a contact manifold, $M$ is a manifold
with boundary of the same dimension as $\tilde{M}$ embedded in
$\tilde{M}$.  Suppose further that there is a neighborhood $U$ in
$\tilde M$ of the boundary $\partial M$ and a free $S^1$ action on $U$
preserving $\alpha$ such that the corresponding moment map $f:U \to
\R$ satisfies
\begin{enumerate}
\item $f\inv (0) = \partial M $ and
\item $f\inv ([0, \infty)) = U \cap M$.
\end{enumerate}
 
Let $M_{\text{cut}} = M/ \sim$, where, for $m\not = m'$, $m\sim m'$ if
and only if
\begin{enumerate}
\item $m, m' \in \partial M$ and 
\item $m = \lambda \cdot m'$ for some $\lambda \in S^1$, 
\end{enumerate}
Then $M_{\text{cut}}$ is a contact manifold, $\partial M/S^1$ is a
contact submanifold of $M_{\text{cut}}$, and $M_{\text{cut}}
\smallsetminus (\partial M/S^1)$ is contactomorphic to
$M\smallsetminus \partial M$.

Moreover if there is an action of a Lie group $G$ on $\tilde{M}$
preserving $M$, $\alpha$ and commuting with the action of $S^1$ on
$U$, then there is an induced action of $G$ on $M_{\text{cut}}$
preserving the induced contact structure.
\end{proposition}

Suppose we are given $t_1, t_2 \in \R$ with $\tan t_1, \tan t_2$
rational, $0\leq t_1 < 2\pi$ and $t_1 < t_2$.  For each $i$ there is a
weight $(m_i,n_i) \in \Z^2$ such that $(\cos t_i, \sin t_i)$ lies on
the ray through $(m_i, n_i)$.

Choose $\epsilon > 0$ sufficiently small so that $f_1 (t) = - n_1 \cos
t + m_1 \sin t$ is non-negative on $[t_1, t_1 + \epsilon)$ and $f_2
(t) = n_2 \cos t - m_2 \sin t$ is non-negative on $(t_2 -\epsilon,
t_2]$ and $t_1 + \epsilon < t_2 - \epsilon$.  Consider $\tilde{M} =\R
\times S^1 \times S^1$ with the contact form $\alpha = \cos t \,
d\theta_1 + \sin t\, d\theta_2$, $(t,
\theta_1, \theta_2)\in \R \times S^1 \times S^1$.  Let 
$M = [t_1, t_2] \times S^1 \times S^1$, 
$U = ((t_1-\epsilon, t_1 + \epsilon) \cup (t_2 -\epsilon, t_2 + \epsilon)) \times S^1 \times S^1$.  Consider $f: U\to \R$ given by 
$f(t, \theta_1, \theta_2) = - n_1 \cos t + m_1 \sin t$ for $t\in
(t_1-\epsilon, t_1 + \epsilon)$ and $f(t, \theta_1, \theta_2) =
n_2\cos t + m_2 \sin t$ for $t\in (t_2-\epsilon, t_2 + \epsilon)$.
The function $f$ is a moment map for a circle action on $U$ generated
on $(t_1-\epsilon, t_1 + \epsilon) \times S^1 \times S^1$ by $-n_1
\frac{\partial}{\partial \theta_1} + m_1 \frac{\partial}{\partial \theta_2}$ 
and on $(t_2-\epsilon, t_2 + \epsilon) \times S^1 \times S^1$ by $n_2
\frac{\partial}{\partial \theta_1} - m_2 \frac{\partial}{\partial \theta_2}$.

Note that the obvious action of $G= S^1 \times S^1$ on $\tilde{M}$
preserves $M$, $\alpha$, $U$ and commutes with the action of $S^1$
defined by $f$.  Therefore we may apply Proposition~\ref{cont-prop2}.
The contact manifold $M_{\text{cut}}$ so obtained with the induced
action of $G$ is the desired manifold $(M, \alpha, \Psi_\alpha)$.
\end{proof}

\subsection{\protect Proof of Theorem~\ref{main-theorem} (3)}

Let $(M, \alpha, \Psi_\alpha)$ be a c.c.c.t. $G$-manifold normalized
so that $\Psi_\alpha (M) \subset S(\fg^*)$, suppose the action of $G$
is free and suppose $\dim M > 3$.  By Corollary~\ref{Cor1} the moment
map $\Psi_\alpha : M\to S(\fg^*)$ is a principal $G$-bundle.  By
Lemma~\ref{BlockA} $(M, \alpha, \Psi_\alpha)$ is locally isomorphic to
the co-sphere bundle $ S^*G= G\times S(\fg^*)$ of the torus $G$ with
the standard contact structure and the obvious action of $G$.  By
Proposition~\ref{prop1*} $(M, \alpha, \Psi_\alpha)$ corresponds to a
class in $H^1(M/G, \cS) = H^1 (S(\fg^*), \cS)$.  By
Corollary~\ref{cor4.4} $H^1 (S(\fg^*), \cS)$ is isomorphic to $H^2
(S(\fg^*), \underline{\Z}_G)$.  On the other hand the cohomology group $H^2
(S(\fg^*), \underline{\Z}_G)$ classifies principal $G$-bundles over the sphere
$S(\fg^*)$.  This suggests (but doesn't yet prove!) that every
principal $G$-bundle over $S(\fg^*)$ has a unique invariant contact
structure.  To prove this we need to trace through identifications.

Recall that principal $G$-bundles (for $G$ a torus) over $B= S(\fg^*)$
are in 1-1 correspondence with classes in the first \v{C}ech
cohomology $H^1 (B, \underline{G})$ where $\underline{G}$ is the sheaf
defined by $\underline{G}(U) = C^\infty (U, G)$, for $U\subset B$
sufficiently small.  Recall also that we have a short exact sequence
of sheaves 
$$ 
0 \to \underline{\Z}_G \to \underline{\fg} \stackrel{\exp}{\to}
\underline{G} \to 1 , 
$$ 
where $\underline{Z}_G (U) = C^\infty (U, \Z_G) \simeq \Z_G$,
$\underline{\fg} (U) = C^\infty (U, \fg)$ and $\exp: \underline{\fg}
\to \underline{G}$ is induced by $\exp : \fg \to G$.
On the other hand, by Proposition~\ref{prop4.3} we have the short
exact sequence 
$$ 
0 \to \underline{\Z}_G \stackrel{j}{\to} \cC
\stackrel{\Lambda}{\to} \cS \to 1.  
$$ 
We claim that there are morphisms $a: \cC \to \underline{\fg}$ and $b:
\cS \to \underline{G}$ so that the diagram
\begin{equation*}
\begin{CD}
 \underline{\Z}_G @>j>> \cC @>\Lambda>> \cS \\
@VV{id}V @VaVV @VbVV \\
\underline{\Z}_G @>>> \underline{\fg} @>>> \underline{G}
\end{CD}
\end{equation*}
commutes and that, moreover, $b$ induces an isomorphism $H^1
(S(\fg^*), \cS) \to H^1 (S(\fg^*), \underline{G})$.

Denote the projection $M \to S(\fg^*)$ by $\pi$.  Fix an open set
$U\subset S(\fg^*)$ small enough  so that $\pi\inv (U) = U \times G$.  By
Lemma~\ref{lemma?} for any $f\in \cC (U)$ the contact vector field
$\Xi _f$ of $f$ is tangent to $G$-orbits, i.e., to the fibers of
$\pi$. Since the vector field $\Xi _f$ is also $G$-invariant, for any
$x\in U$ there is a unique $X(x) \in \fg$ such that for any $m\in \pi\inv (x)$
$$
\Xi _f (m)  = (X(x))_M (m)
$$
(cf.\ proof of Proposition~\ref{prop4.3}). 
We define $a(f)$ to be the map $X:U\to \fg$, $x\mapsto X(x)$.
Since an element $\varphi \in \cS (U)$ is a $G$-equivariant
diffeomorphism of $\pi\inv (U) = U \times G$ into itself, it is
completely determined by $\varphi|_{U\times \{1\}}$. We define $b: \cS
(U) \to \underline{G}$ by $b(\varphi)(x) = \varphi(x, 1)$.  

By definition of $\Lambda : \cC \to \cS$, $\Lambda (f)$ is the time 1
flow of $\Xi_f$.  Thus if $\Xi _f (m) = (X(\pi(m)))_M (m)$ then
$\Lambda (f) (m) = (\exp X(\pi (m)))\cdot m = \exp (a(f) (\pi
(m)))\cdot m$.  Consequently $b (\Lambda (f)) (x) = \exp (a(f)(x))$,
i.e., $ b \circ \Lambda = \exp \circ a$.

Finally the left hand square commutes by definition of $j$ and the
fact that the contact vector field of the function $\langle
\Psi_\alpha , X\rangle$ is $X_M$.  This proves the claim. 

Since $\cC$ and $\underline{\fg}$ are fine sheaves we have
\begin{equation*}
\begin{CD}
H^1 (S(\fg^*), \cS ) @>\cong>> H^2 (S(\fg^*), \underline{\Z}_G  )\\
@V{H^1 (b)}VV     @V{H^2(id)}VV     \\
H^1 (S(\fg^*), \underline{G})  @>\cong>> H^2 (S(\fg^*), \underline{\Z}_G  )\\
\end{CD} .
\end{equation*}
Therefore the map $H^1(b)$ induced by $b$ on the first \v{C}hech
cohomology is an isomorphism.

\subsection{\protect Proof of Theorem~\ref{main-theorem} (4)}

\begin{lemma} 
Let $\fg^*$ be the dual of the Lie algebra of a torus $G$.  Let
$C\subset \fg^*$ be a good polyhedral cone.  There exists a
c.c.c.t. $G$-manifold $(M, \alpha, \Psi_\alpha)$ such that the moment
cone of $\Psi_\alpha$ is $C$.
\end{lemma}
\begin{proof}
Suppose $C =\bigcap \{ \eta \in \fg^* \mid \langle \eta, v_i\rangle
\geq 0\}$ is a good polyhedral cone defined by some subset $\{v_i\}_{i=1}^N$
of the integral lattice $\Z_G$.

As a first step we construct a symplectic cone $(S,\omega)$\footnote{
Recall that a symplectic cone is a symplectic
manifold $(S, \omega)$ with a free proper action $\{\rho_t\}_{t\in \R}$
of the real line such that $\rho_t^* \omega = e^t \omega$ for all $t\in
\R$} with a
symplectic action of $G$ commuting with dilations such that the image
of the corresponding moment map $\Phi_S: S \to \fg^*$ is $C$. The
construction is a slight adaptation of a well-known construction of
Delzant (c.f.\ \cite{D}, \cite{LT}).

Let $\{e_i\}$ denote the standard basis of $\R^N$.  Consider the map
$\varpi: \R ^N \to \fg$ given by $\varpi (\sum a_i e_i) = \sum a_i
v_i$. Since $\varpi (\Z^N) \subset \Z_G$, $\varpi$ induces a map
$\bar{\varpi}: \bbT^N = \R^N/\Z^N \to \fg/\Z_G = G$.  We write $[a]$
for the image of $a =(a_1, \ldots, a_N)\in \R^N$ in $\bbT^N$.  Note
that the kernel $T$ of $\bar{\varpi}$ is 
$$ 
T = \{ [a] \mid \sum a_i v_i \in \Z_G \}.  
$$ 
It is a compact abelian subgroup of $\bbT^N$  with
Lie algebra $\ft = \ker \varpi$.  Note that $T$ need not be connected.

Consider the standard action of $\bbT^N$ on $(\C^N,
\frac{\sqrt{-1}}{2\pi} \sum dz_j \wedge d \bar{z} _j)$:
$$
[a]\cdot (z_1, \ldots, z_N) = (e^{2\pi i a_1}z_1, \ldots e^{2\pi i
 a_N}z_N).
$$  
The corresponding symplectic moment map $\Phi : \C^N \to (\R^N)^*$ is
given by $\Phi (z_1, \ldots, z_N) = \sum |z_j|^2 e_j^*$ where
$\{e_j^*\}$ is the basis dual to $\{e_j\}$.  We claim that the
symplectic quotient $S$ of $\C^N \smallsetminus \{0\}$ by the induced
action of $T$ is the desired manifold $S$, i.e., that $$ S =
(\Phi_T\inv (0) \smallsetminus \{0\})/T $$ where $\Phi_T = j^* \circ
\Phi$ and $j: \ft \hookrightarrow \R^N$ is the inclusion.

We claim first that $\Phi_T \inv (0) = \Phi \inv (\varphi (C))$.
Indeed, since $0 \to \fg^* \stackrel{\varpi^*}{\to} (\R^N)^*
\stackrel{j^*}{\to } \ft^* \to 0$ is exact, $(j^*)\inv (0) = \varpi ^*
(\fg^* )$ and hence $\Phi_T\inv (0) = \Phi\inv ((j^*)\inv (0)\cap \Phi
(\C^N)) = \Phi\inv (\varpi^* (\fg^*) \cap \Phi (\C^N))$.  Now 
\begin{equation*}
\begin{split}
\varpi^* (\fg^*) \cap \Phi (\C^N) 
&= \{ \varpi^* (\eta) \mid \eta \in \fg^*\,\, \text{ and }\,\, 
\langle  \varpi^* (\eta), e_i\rangle \geq 0 \text{ for all } i\}\\
&= \{ \varpi^* (\eta) \mid \eta \in \fg^*\,\, \text{ and }\,\, 0\leq
\langle \eta, \varpi (e_i) \rangle = \langle \eta, v_i\rangle \, 
\text{ for all } i\}\\
&= \{ \varpi^* (\eta) \mid \eta \in C\}.
\end{split}
\end{equation*}
Thus $\Phi_T\inv (0) = \Phi \inv (\varpi^* (C))$.

Next we claim that for any $0\not = z= (z_1, \ldots, z_N) \in \Phi_T
\inv (0)$ the isotropy group $T_z$ is trivial.  It would then follow
that $S = (\Phi_T\inv (0) \smallsetminus \{0\})/T$ is a smooth
symplectic cone: the action of $\R$ on $S$ is induced by the action of
$\R$ on $\C^N$ given by $t \cdot z = e^t z$.  Now $T_z = T \cap
(\bbT^N)_z = \{ [a]\in \R^N/\Z^N \mid \sum a_i v_i \in \Z_G\} \cap
\{[a]\in \R^N/\Z^N \mid a_i \in \Z \, \text{ for all $i$ with } z_i
\not = 0\} = \{[a]\in \R^N/\Z^N \mid \sum _{j\in J_z} a_j v_j \in \Z_G
\text{ and } a_j \in \Z \text{ for all } j\not\in J_z\}$ where $J_z =
\{j\in \{1, \ldots, N\} \mid z_j = 0\}$.

On the other hand, $z\in \Phi_T \inv (0)$ if and only if there is a
(unique) $\eta \in C$ such that $\Phi (z) = \varpi^* (\eta)$.  Hence
$z_j= 0$ if and only if $0= |z_j|^2 = \langle \Phi (z), e_j\rangle =
\langle \varpi^* (\eta), e_j \rangle = \langle \eta, v_j\rangle$.  Since $C$ 
is a good cone, for any fixed vector $\eta \in C$ the set $\{v_j \mid
\langle \eta , v_j\rangle =0\} = \{ v_j \mid j\in J_z\}$ is a $\Z$
basis of $\{ \sum _{j\in J_z} a_j v_j \mid a_j \in \R\} \cap \Z_G$.
Hence $\sum _{j\in J_z} a_j v_j \in \Z_G$ implies that $a_j \in \Z$
for all $j\in J_z$. Therefore $T_z = \{ [a] \in \R^N/\Z^N \mid a\in
\Z^N\}$, i.e., $T_z$ is trivial.

Finally note that the image of $\Phi_T \inv (0)$ under $\Phi$ is
precisely $\varpi ^* (C)$.  Hence the image of the reduced space
$S = \Phi_T\inv (0)/T$ under the induced $\bbT^N/T = G$ moment map
$\tilde{\Phi} :\Phi_T\inv (0)/T \to \ft^\circ = \varpi^* (\fg^*)$ is
$\varpi^*(C) \simeq C$.

Since the sphere $S^{2N-1} = \{z\in \C^N \mid ||z||^2 = 1\}$ is a 
$\bbT^N$-invariant hypersurface of contact type in $\C^N$, and since 
the action of $\R$ on $\C^N$ commutes with the action of $\bbT^N$
$$ 
M := (\Phi_T\inv (0) \cap S^{2N-1})/T 
$$ 
is a $\bbT^N/T = G$-invariant hypersurface of contact
type in the quotient $\Phi_T\inv (0)/T$.  Moreover $\tilde{\Phi}|_M: M
\to \ft^\circ = \varpi^* (\fg^*)$ is the corresponding contact moment
map, and its moment cone is precisely $\varpi^* (C) \simeq C$.
\end{proof} 

\begin{lemma}
Suppose $(M, \alpha, \Psi_\alpha :M \to \fg^*)$ is a
c.c.c.t. $G$-manifold, $\dim M > 3$ and the action of $G$ is not free.
Then the moment cone $C(\Psi)$ is a good rational polyhedral
cone.
\end{lemma}

\begin{proof}
We first introduce some notation and a simple fact.  Let $C\subset
\fg^*$ be a cone and $F\subset C$ be a face of $C$.  Let
$\text{span}_\R F$ denote the linear subspace of $\fg^*$ spanned by
the vectors in $F$.  Let $\pi_F : \fg^* \to \fg^* /\text{span}_\R F $
denote the projection.  For any point $q$ in the interior of $F$ there
is an open neighborhood $\W$ of $q$ in $\fg^*$ such that
$$
C \cap \W = \pi_F \inv (\pi_F (C)) \cap \W.
$$
Note that the cone $\pi_F \inv (\pi_F (C))$ is isomorphic to $\pi_F
(C) \times\text{span}_\R F$.

Now suppose that $C$ is the moment cone $C(\Psi)$ of
$\Psi_\alpha$ and that $x$ is a point in $M$. By Lemma~\ref{lemma-*}
for any $G$ invariant neighborhood $U$ of $G\cdot x$ there is an open
subset $\tilde{W}$ of the sphere $S(\fg^*)$ such that $\Psi_\alpha (U)
= \tilde{W} \cap \Psi_\alpha (M)$.  Let $W = \R^+ \tilde{W} \cup
\{0\}$ be the cone on $\tilde{W}$.  Then
$$
C(\Psi_\alpha|_U) = C(\Psi_\alpha) \cap W.
$$
By the local normal form theorem, Theorem~\ref{BlockF},
$$
C(\Psi_\alpha|_U) = W \cap \left( 
\R^+ \left(\Psi_\alpha (x) +j ((\fk/\fg_x)^*) + i (\Phi_V (V))\right) 
\cup \{0\} \right),
$$ 
where as usual $\fk$ is the characteristic subalgebra, $\Phi_V :V
\to \fg^*_x$ is the moment map for the slice representation etc..
Note that 
$$
 W \cap  
\R^+ (\Psi_\alpha (x) +j ((\fk/\fg_x)^*))  = W \cap \fg_x^\circ
$$
for any sufficiently ``small'' open cone $W$.  Thus 
$$ 
C(\Psi_\alpha|_U) = W
\cap (\fg_x^\circ + i (\Phi_V (V))) = W \cap \pi_x \inv (\Phi_V (V)), 
$$ 
where $\pi_x: \fg^* \to \fg_x^*$ is the natural projection. It
follows that if $F$ is the face of $C(\Psi_\alpha)$ containing
$\Psi_\alpha (x)$ in its interior, then $\text{span}_\R F =
\fg_x^\circ$ and $\pi_F (C(\Psi_\alpha))$ is isomorphic to $\Phi_V(V)$
(once we identify $\fg^*/\fg_x^\circ$ with $\fg_x^*$).

Moreover, if we represent $C$ as 
$C =\bigcap \{ \eta \in \fg^* \mid \langle \eta,
v_i\rangle \geq 0\}$ for some minimal  set $\{v_i\} \subset \Z_G$ 
and consisting of primitive vectors, then $\fg_x^\circ = \text{span}_\R F
= \bigcap _{j=1}^k \{ \eta \in \fg^* \mid \langle \eta, v_{i_j}
\rangle = 0\}$ for some subset $\{v_{i_1}, \ldots , v_{i_k}\} \subset
\Z_G$, and
$$
W \cap C(\Psi_\alpha)= W\cap \left(\pi_x\inv (\Phi_V(V))\right) = 
W \cap \bigcap _{j=1}^k \{ \eta \in \fg^* \mid \langle \eta, v_{i_j}
\rangle \geq 0\}.
$$ 
Hence $\pi_x\inv (\Phi_V (V)) = \bigcap _{j=1}^k \{ \eta \in \fg^*
\mid \langle \eta, v_{i_j}\rangle \geq 0\}$.  Since 
$\Z_{G_x} = \fg_x \cap \Z_G$, $\{v_{i_j}\}$ is a subset of
$(\fg_x^\circ)^\circ \cap \Z_G = \fg_x \cap \Z_G = \Z_{G_x}$.  For any
$v\in \fg_x$ we have $\pi_x (\{ \eta\in \fg^* \mid \langle \eta, v\rangle \geq
0\}) =\{\eta\in \fg^*_x \mid \langle \eta, v\rangle \geq 0\})$.  
Therefore 
$$
\Phi_V (V) = \pi_x \left( \bigcap _{j=1}^k \{ \eta \in \fg^* 
\mid \langle \eta, v_{i_j}\rangle \geq 0\}\right) =
\bigcap _{j=1}^k \{ \eta \in \fg^*_x \mid \langle \eta, v_{i_j}
\rangle \geq 0\}
$$ 
By minimality of $\{v_i\}$ the set $\{v_{i_j}\}$ is the minimal set
with this property.
On the other hand, by Lemma~\ref{lemma.conn+slice},
$$
\Phi_V (V) = \{ \sum _{j=1}^k a_j \nu_j \mid a_j \geq 0\}
$$ for some basis $\{\nu_j\}$ of $\Z_{G_x}^*$.  Therefore the set
$\{v_{i_j}\}$ is a basis of the lattice $\Z_{G_x}$.
\end{proof}

\section{Uniqueness of toric integrable actions  on $T^* \bbT ^n$ and 
on $T^* S^2$}

Having proved the main classification theorem,
Theorem~\ref{main-theorem}, we now in position to prove
Theorems~\ref{thm1.1} and \ref{thm1.2}.

\begin{proof}[Proof of Theorem~\ref{thm1.1}]
Let $(M, \alpha, \Psi_\alpha : M \to \fg^*)$ be a contact toric
$G$-manifold ($G = \bbT^n$) such that $(M, \alpha)$ is contactomorphic
to the co-sphere bundle $S^* \bbT^n$ of the $n$-torus $\bbT^n$ with
the standard contact structure.  We will argue that $(M, \alpha,
\Psi_\alpha)$ is unique as a contact toric manifold.

It was shown in \cite{LS} that if $M = S^*\bbT^n$, then the action of
$G$ is necessarily free.  The argument roughly goes as follows (see
\cite{LS} for details).  Suppose the action of $G$ is not free.
Consider first the case of $\dim M = 3$.  Then $M$ is a lens space
(cf.\ Theorem~\ref{main-theorem}~(2)), hence cannot be $S^* \bbT^2 =
\bbT^3$.  

Next consider the case of $\dim M > 3$.  Then the moment cone $C$ of
$(M, \alpha, \Psi_\alpha)$ is a good polyhedral cone determining $(M,
\alpha, \Psi_\alpha)$ uniquely (cf.\ Theorem~\ref{main-theorem}~(4)).
If the maximal linear subspace of $C$ has dimension $k>0$ then $C$ is
isomorphic to the moment cone of $M' = \bbT^k \times S^{2n -1 -k}$,
where $M'$ gets its contact toric structure as a hypersurface of
contact type $\{ (q, p, z) \in \bbT^k \times (\R ^k)^* \times \C^{n-k}
\mid \, ||p||^2 + ||z||^4 = 1\}$ in $T^*\bbT^k \times \C^{n-k}$.  Consequently
 $M$ is homeomorphic to $M' = \bbT^k \times S^{2n -1 -k} 
\not = \bbT^n \times S^{n-1} = S^* \bbT^n$.

Finally if the dimension of the maximal linear subspace of $C$ is
zero, i.e., if $C$ is a proper cone, then by a theorem of Boyer and
Galicki \cite{BG} $M$ has a locally free circle action so that the
quotient $M/S^1$ is a compact connected symplectic toric orbifold.
The real odd-dimensional cohomology of a compact symplectic
toric orbifold is zero.  Consequently $\dim H^1 (M, \R) \leq 1$.  Hence $M
\not = S^* T^n$, $n>1$.

We conclude that if $(M, \alpha, \Psi_\alpha : M \to \fg^*)$ is a
contact toric $G$-manifold ($G = \bbT^n$) and $M = S^* \bbT^n$ then
the action of $G$ is necessarily free.  We  argue next that it
is unique.

Suppose $\dim M = 3 $ and the action of $G= \bbT^2$ is free.  By the
classification theorem $(M, \alpha) = (\bbT^3, \alpha _k = \cos kt
\,d\theta _1 + \sin kt \, d\theta_2)$, $k = 1, 2, \ldots$.  By a
theorem of Giroux \cite{Gi}, $(\bbT^3, \alpha_k)$ and $(\bbT^3,
\alpha_l)$ are distinct as contact manifolds for $k\not =l$.  Since 
$\alpha_1$ is the standard contact structure on $S^* \bbT^2$, it
follows that there is only one contact toric manifold contactomorphic
to $(S^*\bbT^2, \alpha_1)$.  In other words there is only one
$\bbT^2$-action on $S^*\bbT^2$ making it a contact toric manifold.

Suppose next that $\dim M >3$ and the action of $G$ is free. By
Theorem~\ref{main-theorem}~(4), $M$ is a principal $G$-bundle over the
sphere $S^{n-1}$, $n = \dim G$, and each such principal $G$-bundle
carries only one $G$-invariant contact structure.  Now, principal
$\bbT^n$ bundles over $S^{n-1}$ are in one-to-one correspondence with
elements of $H^2 (S^{n-1}, \Z^n)$ which is 0 unless $n-1 = 2$, in
which case it's $\Z^3$.  Note however that no nontrivial $\bbT^3$
bundle over $S^2$ is homeomorphic to $S^2 \times \bbT^3$.  We conclude
that if $(M, \alpha, \Psi_\alpha)$ is a contact toric $G$-manifold
such that the action of $G$ is free, $M = S^* G$ and $\dim M > 3$
then $(M, \alpha, \Psi_\alpha)$ is a unique contact toric $G$-manifold
with such properties.  In other words there is only one
$\bbT^n$-action on $S^* \bbT^n$ making it a contact toric manifold.
\end{proof}

\begin{proof}[Proof of Theorem~\ref{thm1.2}]

Suppose $\tau_1$, $\tau_2$ are two effective actions of the torus $G=
\bbT^2$ on $M= S^* S^2 = \R P^3$ preserving the standard contact
structure.  Let $\Psi_1, \Psi_2 :M \to S^1 \subset\fg^*$ be the moment
maps for the actions corresponding to a normalized contact form $\alpha_0$
defining the standard contact structure.  We will argue that the
images $\Psi_i (M)$ are arcs in $S^1$ of length less than $\pi$ (hence
$\Psi_i$ are one-to-one).  Moreover we'll show that there is an
element $A \in \text{SL}(\fg^*)$ preserving the weight lattice
$\Z_G^*$ such that $A (\Psi_1 (M)) = \Psi_2 (M)$.  It would then
follow that the action $\tau_1$ composed with the isomorphism of $G$
defined by $A$ is $\tau_2$.

Note that since $\R P^3 \not = \bbT^3$, the actions $\tau_i$ are
necessarily not free (cf.\ Theorem~\ref{main-theorem} and the proof of
Theorem~\ref{thm1.1} above).  Now consider one of the two actions, say
$\tau_1$.  By Lemma~\ref{lemma6.1} the action is free except at two
orbits $G\cdot x_1$ and $G\cdot x_2$.  The isotropy groups $K_i$ of
$x_i$ are circles, and the images $\Psi_1 (x_i)$ are of the form
$\frac{\mu_i}{||\mu_i||}$ where $\mu_i \in \Z_G^*$ are primitive
weights with $\ker \mu_i$ being the Lie algebra of the circle $K_i$.

It follows from the proof of Lemma~\ref{lemma6.1} that the contact
toric manifold $(M, \alpha_0, \Psi_1 :M\to \fg^*)$ can be obtained by
cutting $(\tilde{B}, \tilde{\alpha}) = (\bbT^2 \times \R, \cos t \,
d\theta_1 + \sin t \, d\theta_2)$ using $B = \bbT^2 \times [t_1,
t_2]$, where $(\cos t_i, \sin t_i) = \frac{\mu_i}{||\mu_i||}$ (we
identify $\fg^*$ with $\R^2$ and $\Z_G^* $ with $\Z^2$).  Hence as a
topological space $M = B/\sim$ where $(g, t_1)\sim (ag, t_1)$ for all
$a\in K_1$ and $(g, t_2)\sim (ag, t_2)$ for all $a\in K_2$.  Note that
$t_2 - t_1 \not = \pi n$ , $n=1, 2, \ldots$, for otherwise $\mu_1 =
\pm \mu_2$ and then $B/\sim$ is $S^2 \times S^1 \not = \R P^3$.

Next observe that since the standard contact structure on $S^* S^2 =
M$ is tight, we must have $t_2 - t_1 < \pi$ (cf.\
\cite{L-IJM}). Indeed, if $t_1 = 0$ and $t_2 > \pi$, the image of the
cylinder $\{(1, e^{i\theta}, t) \mid 0\leq t \leq t_2,\, \theta \in
\R\} \subset B \subset \bbT^2 \times \R$ in $M = B/\sim$ is an
overtwisted disk.  One can write a similar formula for an overtwisted
disk if $t_1 >0$.  We conclude that the image $\Psi_1 (M)$ in $S^1$ is
an arc of length less than $\pi$.  Consequently the fibers of $\Psi_1
:M \to S^1$ are connected.

We next argue that the weights $\mu_1$, $\mu_2$ which span the edges
of the moment cone $C(\Psi_1)$ span a sublattice of the weight lattice
$\Z_G^*$ of index two.

\begin{lemma}
 Let $G= \bbT^2$ and let 
$K_1, K_2 \subset G$ be two closed subgroups isomorphic to $S^1$.  Let
$M$ be the topological space $(G\times [0, 1] /\sim$ where $(0, g)
\sim (0, ag)$ for all $g\in G$, $a \in K_1$ and $(1, g) \sim (1, ag) $
for all $g\in G$ and $a\in K_2$.  In other words $M$ is obtained from
the manifold with boundary $G \times [0,1]$ by collapsing circles in
the two components of the boundary by the respective actions of two
 circle subgroups.
Let $\mu_1, \mu_2 \in Z_G^*$ be the two primitive weights determined by
$K_1$ and $K_2$ respectively, i.e., the kernel of the character
defined by $\mu_i$ is $K_i$.

Then $H^1 (M, \Z) \simeq \{ (n_1, n_2) \in \Z^2 \mid n_1 \mu_1 + n_2
\mu_2 = 0\}$ and $H^2 (M, \Z) \simeq \Z_G^* /(\Z\mu_1 + \Z \mu_2)$.
\end{lemma}

\begin{proof}
Recall that $H^1 (G, \Z)$ is isomorphic to the weight lattice $\Z_G^*$
and that the isomorphism is given as follows.  A weight $\nu \in
\Z_G^*$ defines a character $\chi_\nu : G \to S^1$ by $\chi_\nu ( \exp
(X)) = e^{2\pi i \nu (X)}$; the class $\chi_\nu ^* [d\theta]$ is the
element in $H^1 (G, \Z)$ corresponding to $\nu$.  Here $d\theta$ is
the obvious 1-form on $S^1$.

Consequently if $G= \bbT^2$ and $K_j \subset G$ is a circle subgroup,
then $\pi_j : G \to G/K_j \simeq S^1$ is a character and hence the
weight $\mu_j = (d \pi_j)_1 $ defines an element of $H^1 (G, \Z)$.
Thus if we identify $H^1 (G/K_j, \Z)$ with $\Z$ and $H^1 (G, \Z)$ with
$\Z_G^*$, then the map $H^1 (G/K_j, \Z) \to H^1 (G, \Z)$ becomes the
map $\Z\ni n\mapsto n\mu_j \in \Z_G^*$.

The sets $U = (G \times [0, 2/3) )/\sim $ and $V = (G \times(1/3, 1]
)/\sim$ are two open subsets of $M$.  We have $M = U \cup V$, $U\cap V
= G \times(1/3, 2/3)$ is homotopy equivalent to $G$, $U$ is homotopy
equivalent to $G/K_1$, $V$ is homotopy equivalent to $G/K_2$ and the
inclusion maps $U\cap V \hookrightarrow U$, $U\cap V \hookrightarrow
V$ are homotopy equivalent to projections $\pi_1: G \to G/K_1$, $\pi_2
:G \to G/K_2$ respectively.  Hence under the above identifications of
$H^1 (U)$ and $H^1(V)$ with $\Z$, the inclusions $U\cap V \to U$,
$U\cap V
\to V$ induce the maps $\Z\ni n\mapsto n\mu_j \in \Z_G^*$, $j=1, 2$,
respectively.

We now apply the Mayer-Vietoris sequence to compute the integral
cohomology of $M$.  We start with $0 \to H^0(M) \to H^0(U) \oplus H^0
(V) \to H^0 (G)\stackrel{\delta}{ \to} H^1(M) \to H^1(U) \oplus H^1(V)
\to H^1 (G)
\stackrel{\delta}{ \to} H^2 (M)\to H^2(U) \oplus H^2(V) \to H^2 (G)
\stackrel{\delta}{ \to}H^3 (M) \to 0$.  Clearly the map $H^0(U) \oplus
H^0 (V) \to H^0 (G) $ is onto.  Given the identifications above the
map $\varphi :H^1(U) \oplus H^1(V) \to H^1 (G)$ becomes $\Z \oplus \Z\ni
(n, m) \mapsto n\mu_1 + m \mu_2 \in \Z_G^*$.   We therefore have
$0 \to H^1 (M) \to \Z \oplus \Z \stackrel{\varphi}{\to } \Z_G^*
\stackrel{\delta}{ \to} H^2 (M) \to 0\oplus 0 \to H^2 (G) 
\stackrel{\delta}{ \to} H^3 (M) \to 0$ and the lemma follows.
\end{proof}

Since $M= \R P^3$, $H^2 (M, \Z) = \Z/2$.  Hence $\Z/2 = \Z_G^*
/(\Z\mu_1 + \Z \mu_2)$.  Consequently, since $\mu_1, \mu_2$ are
primitive, the parallelogram 
$$ 
\{ a_1 \mu_1 + a_2 \mu_2 \mid 0\leq a_1, a_2 \leq 1\} 
$$ 
contains exactly five point of $\Z^*_G$: four vertices plus the point
$\mu = \frac{1}{2} (\mu_1 + \mu_2)$ in its interior.  Hence $\{\mu_1,
\mu\}$ is a basis of $\Z_G^*$.  Of course $\mu_2 = 2 \mu - \mu_1$.

By the same argument the image $\Psi_2 (M) $ is an arc in $S^1$ of
length less than $\pi$ with endpoints $\frac{\nu_1}{||\nu_1||}$,
$\frac{\nu_2}{||\nu_2||}$where $\nu_1, \nu_2 \in \Z_G^*$ are primitive
weights.  Moreover $\{\nu_1, \nu = \frac{1}{2} (\nu_1 + \nu_2)\}$ is a
basis of $\Z_G^*$.  The linear map $A: \fg^* \to \fg^*$ defined by
$A\mu_1 = \nu_1$, $A\mu = \nu$ is the desired map.
\end{proof}


\begin{thebibliography}{WWWW}
\bibitem[A]{Albert} C. Albert, 
Le th\'eor\`eme de r\'eduction de Marsden-Weinstein en g\'eom\'etrie
cosymplectique et de contact, {\em J.\ Geom.\ Phys.} {\bf 6} (1989),
no.~4, 627--649.

\bibitem[Bn]{B} A. Banyaga, 
The geometry surrounding the Arnold-Liouville theorem in {\em Advances
in geometry}, edited by Jean-Luc Brylinski, Ranee Brylinski, Victor
Nistor, Boris Tsygan and Ping Xu. Progress in Mathematics, 172.
Birkhäuser Boston, Inc., Boston, MA, 1999. xii+399 pp.  ISBN
0-8176-4044-4

\bibitem[BnM1]{BM1}A. Banyaga and P. Molino
G\'eom\'etrie des formes de contact compl\`etement int\'egrables de
type toriques in {\em   S\'eminaire Gaston Darboux de G\'eom\'etrie et
Topologie Diff\'erentielle}, 1991--1992 (Montpellier), 1--25,
Univ. Montpellier II, Montpellier, 1993.

\bibitem[BnM2]{BM2} A. Banyaga and P. Molino, 
Complete integrability in contact geometry, Penn State preprint PM 197, 1996.

\bibitem[Bt]{Bates} L. Bates,  Examples for obstructions to action-angle 
coordinates, {\em  Proc.\ Roy.\ Soc.\
Edinburgh Sect.\ A} {\bf 110} (1988), no. 1-2, 27--30.

\bibitem[BoM]{Bou-M} M. Boucetta and P. Molino,
G\'eom\'etrie globale des syst\`emes hamiltoniens compl\`etement
int\'egrables: fibrations lagrangiennes singuli\`eres et coordonn\'ees
action-angle \`a singularit\'es, {\em
C.\ R.\ Acad.\ Sci.\ Paris S\'er.\ I Math.} {\bf 308} (1989), no. 13, 421--424.

\bibitem[BG]{BG} C. P. Boyer and K. Galicki, A note on toric contact geometry,
 {\em J. of Geom.\ and Phys.} {\bf 35} (2000) 288--298;
{\tt math.DG/9907043v2}.



\bibitem[D]{D} T. Delzant,  Hamiltoniens p\'eriodiques et images convexes de 
l'application moment, {\em Bull.\ Soc.\ Math.\ France} {\bf 116}
(1988), no. 3, 315--339.

\bibitem[Ge]{Geiges} H. Geiges, Constructions of contact manifolds,
{\em Math.\ Proc.\ Cambridge Philos.\ Soc.} {\bf 121} (1997), no. 3,
455--464.

\bibitem[Gi]{Gi} E. Giroux, Une structure de contact, m\^{e}me tendue, est 
plus ou moins tordue, {\em  
Ann.\ Sci.\ \'Ecole Norm.\ Sup.\ (4)} {\bf 27} (1994), no. 6, 697--705.

\bibitem[GS]{GS} V. Guillemin and S. Sternberg, 
{\em Symplectic techniques in physics},  Cambridge University Press,
 Cambridge -- New York, 1984. xi+468 pp. ISBN: 0-521-24866-3

\bibitem[HS]{HS} A. Haefliger and E. Salem, Action of tori on orbifolds, 
{\em Ann.\ Global Anal.\ Geom.} {\bf 9} (1991), 37--59.
 
\bibitem[K]{Ko} J. L. Koszul, Sur certain groupes de transformations de Lie, 
dans {\em Colloque de G\'eom\'etrie Diff\'erentielle, } Colloque du CNRS
{\bf 71} (1953), 137--141.

\bibitem[L1]{L-convex} E. Lerman,
 A convexity theorem for torus actions on contact manifolds, {\em
 Ill. J.\ Math}, to appear \\ {\tt
 http://xxx.lanl.gov/abs/math.SG/0012017}.

\bibitem[L2]{L-IJM} E. Lerman,  Contact cuts, {\em Israel J.\ Math}, {\bf 124} 
(2001), 77--92.\\   {\tt http://xxx.lanl.gov/abs/math.SG/0002041}.\hfill

\bibitem[LS]{LS} E. Lerman and N. Shirokova, Toric integrable geodesic flows,\\
{\tt http://xxx.lanl.gov/abs/math.DG/0011139}.

\bibitem[LT]{LT} E. Lerman and S. Tolman,  Symplectic toric orbifolds, 
{\em Trans.\ A.M.S.} {\bf 349} (1997), 4201--4230.

\bibitem[Lu]{Lutz} R. Lutz,  Sur la g\'eom\'etrie des structures de contact
 invariantes, {\em Ann. Inst. Fourier (Grenoble)}, {\bf 29} (1979),
 no. 1, xvii, 283--306.



\bibitem[TZ]{TZ} J. Toth and S. Zelditch, Riemannian manifolds with
uniformly bounded eigenfunctions, {\em Duke Math Journal}, to
appear. See {\tt http://xxx.lanl.gov/abs/math-ph/0002038}.
 
\end{thebibliography}
\end{document}